\documentclass[twoside]{article}

\usepackage{aistats2019}

\usepackage[round]{natbib}

\renewcommand{\cite}[1]{\citep{#1}}

\usepackage{amsmath,amsbsy,amsgen,amscd,amssymb,amsthm,amsfonts,stmaryrd}
\usepackage{bm}
\usepackage{microtype}
\usepackage{graphicx}

\newtheorem{theorem}{Theorem}

\theoremstyle{definition}

\numberwithin{theorem}{section}

\newcommand{\R}{\mathbb{R}}
\newcommand{\trace}{\operatorname{tr}}
\newcommand{\dist}{\operatorname{dist}}

\newcommand{\norm}[1]{\Vert #1 \Vert}
\newcommand{\ip}[2]{\big\langle #1, \, #2 \big\rangle}
\newcommand{\prox}[1]{\mathrm{prox}_{#1}}
\newcommand{\proj}[1]{\mathrm{proj}_{#1}}
\newcommand{\lmo}[1]{\mathrm{lmo}_{#1}}

\usepackage{algorithm,algorithmicx,algpseudocode,setspace}  	
\usepackage[nameinlink]{cleveref}
\usepackage{enumitem}

\begin{document}

\runningtitle{A conditional gradient-based augmented Lagrangian framework}
\runningauthor{Alp Yurtsever, Olivier Fercoq, Volkan Cevher}

\twocolumn[

\aistatstitle{A conditional gradient-based augmented Lagrangian framework}

\aistatsauthor{ Alp Yurtsever$^\dagger$ ~~~ Olivier Fercoq$^\ddagger$ ~~~  Volkan Cevher$^\dagger$ \vspace{0.5em}}

\aistatsaddress{ 
$^\dagger$LIONS, Ecole Polytechnique F\'ed\'erale de Lausanne, Switzerland \\
$^\ddagger$LTCI, T\'el\'ecom ParisTech, Universit\'e Paris-Saclay, France
 } ]

\begin{abstract}
\vspace{-1em}
This paper considers a generic convex minimization template with affine constraints over a compact domain, which covers key semidefinite programming applications. The existing conditional gradient methods either do not apply to our template or are too slow in practice. To this end, we propose a new conditional gradient method, based on a unified treatment of smoothing and augmented Lagrangian frameworks. The proposed method maintains favorable properties of the classical conditional gradient method, such as cheap \textit{linear minimization oracle} calls and \textit{sparse representation} of the decision variable. We prove $\mathcal{O}(1/\sqrt{k})$ convergence rate of our method in the objective residual and the feasibility gap. This rate is essentially the same as the state of the art CG-type methods for our problem template, but the proposed method is significantly superior to existing methods in various semidefinite programming applications.
\end{abstract}

\section{INTRODUCTION}
\vspace{-0.5em}

In this paper we focus on the following constrained convex minimization template with affine constraints:
\begin{equation}
\tag{$P$}
\label{eqn:main-template}
\begin{aligned}
&\underset{x}{\text{minimize}} 
& & f(x) + g(B x) \\
&\text{subject to} 
& & x \in \mathcal{X} ~~ \& ~~ A x \in \mathcal{K}
\end{aligned}
\end{equation}
where ${x}$ is the decision variable that live on the convex and compact optimization domain $\mathcal{X} \subseteq \R^n$; $f: \mathcal{X} \to \R$ is a convex differentiable function with $L_f$-Lipschitz continuous gradient; ${A}:\mathcal{X}\to\R^p$ and ${B}:\mathcal{X}\to\R^q$ are known linear maps; $g:\R^q\to\R$ is a convex function which can be non-smooth but we assume that it is $L_g$-Lipchitz continuous; and  $\mathcal{K} \subseteq \R^p$ is a convex set. 

Conditional gradient method (CGM, a.k.a.\ Frank-Wolfe algorithm) has established itself as a scalable method for solving convex optimization problems over structured domains, thanks to its cheaper oracle compared to the projected and proximal gradient methods. This classical method is originated by \citet{FrankWolfe1956}, and its resurgence in machine learning follows \citet{Hazan2012} and \citet{Jaggi2013}. Despite its favorable properties, classical CGM has restrictive assumptions on the problem template such as smoothness of the objective function, and extension of CG-type methods for broader templates is an active research area (cf. \Cref{sec:related-work} for some recent advancements).

\eqref{eqn:main-template} is significantly broader in applications in comparison with the classical CGM template. Non-smooth term $g(Bx)$ not only lets us use regularization to promote additional structures, it can also be used as a non-smooth loss function which generally provides more robustness than smooth functions. Moreover, \eqref{eqn:main-template} has an affine inclusion constraint $Ax \in \mathcal{K}$ and covers standard semidefinite programming (SDP)  with trace constraint in particular. Hence, a large number of problems in machine learning, signal processing, and computer science can be cast within our template, from unsupervised clustering \cite{Peng2007} to generalized eigenvector problems \cite{Boumal2018}, and from maximum cut \cite{Goemans1995} to phase-retrieval \cite{Candes2011}. We refer to Section~5 from \cite{Yurtsever2018} for a detailed discussion on special instances and applications of \eqref{eqn:main-template}.

Affine constraints pose substantial difficulty for first order methods, hence primal-dual methods are typically preferred for solving \eqref{eqn:main-template} in large scale. Among the primal-dual approaches, augmented Lagrangian provides a powerful framework for deriving fast methods. However, these methods rely on proximal-oracles constrained with $\mathcal{X}$.
Unfortunately, this proximal-oracle does not scale well with the problem dimensions and becomes a computational burden for an important part of the applications, in particular for SDP's. 

To this end, we develop conditional gradient augmented Lagrangian framework (CGAL) to exploit more scalable linear minimization oracles ($\lmo{}$):
\begin{equation*}
\lmo{\mathcal{X}}(v) = \arg\min_{x\in\mathcal{X}} \ip{x}{v}.
\end{equation*}
For instance, $\lmo{}$ possesses a rank-$1$ solution that can be efficiently approximated via power method or Lanczos algorithm for matrix factorization problems and SDPs, while other first-order oracles require full dimensional eigen or singular value decompositions.

CGAL can be viewed as a natural extension of the recent method in \cite{Yurtsever2018} from quadratic penalty to an augmented Lagrangian formulation, and it especially focuses on improving the empirical performance. We prove that CGAL converges with $\mathcal{O}(1/\sqrt{k})$ rate both in the objective residual and the feasibility gap. The simplicity of our analysis also enables us to identify adaptive bounds and propose explicit and implementable the dual step-size rules that retain the theoretical convergence rates, while significantly enhancing the practical performance. Our numerical evidence demonstrates superior performance. 

The rest of this paper is organized as follows. We review the notions of smoothing, quadratic penalty and augmented Lagrangian methods in \Cref{sec:prelimineries}. Then, in \Cref{sec:algorithm} we introduce CGAL and the main convergence theorem. We provide detailed discussion and comparison against the existing related work in \Cref{sec:related-work}. Finally, \Cref{sec:experiments} compiles the empirical evidence supporting the advantages of our framework, and \Cref{sec:conclusions} draws the conclusions. Technical details are deferred to the supplementary material. 

\paragraph{Notation.}
We use lowercase letters for vectors (or matrices when considering vector space of matrices), uppercase letters for linear maps, and calligraphic letters for sets. We denote the Euclidean inner product by $\ip{\cdot}{\cdot}$, and the Euclidean norm by $\norm{\cdot}$. We denote the adjoint of a linear map $A$ by $A^\top$. For a set $\mathcal{K}$, its indicator function $\iota_{\mathcal{K}}:\R^q \to \R \cup \{ +\infty\}$ is defined as
\begin{equation*}
\iota_{\mathcal{K}}(z) = {\footnotesize 
\begin{cases}
0 & \text{if $z \in \mathcal{K}$}\\
+\infty & \text{otherwise.}
\end{cases}
}
\end{equation*}

\section{PRELIMINERIES}\vspace{-3mm}
\label{sec:prelimineries}

Our algorithmic design is based on the unified treatment of smoothing, quadratic penalty and augmented Lagrangian frameworks. We review these notions and explain their similarities in this section.

\subsection{Nesterov Smoothing}
 
In the seminal work, \citet{Nesterov2005} introduces a technique for solving some structured non-smooth optimization problems with efficiency estimates $\mathcal{O}(1/\epsilon)$, which is much better than the theoretical lower bound $\mathcal{O}(1/\epsilon^2)$. This technique is known as Nesterov smoothing, and it is widely used in efficient primal-dual methods (e.g., \cite{Nesterov2005eg},\cite{TranDinh2017}).

Nesterov exploits an important class of non-smooth functions $\psi(x)$ that can be written in the max-form: 
\begin{equation*}
\psi(x) = \hat{\psi}(x) +  \max_{u \in \mathcal{U}} \left\{  \ip{Bx}{u} - \hat{\phi}(u) \right\},
\end{equation*}
for some convex and compact set $\mathcal{U} \subset \R^q$, and some continuous convex functions $\hat{\psi}:\mathcal{X} \to \R$ and $\hat{\phi}:\mathcal{U} \to \R$.

Let us consider a prox-function $\delta(u)$ of $\mathcal{U}$, i.e., a strongly convex continuous function on $\mathcal{U}$. Define the center point of this prox-function as
\begin{equation*}
\dot{u} = \arg\min_{u \in \mathcal{U}} \delta(u).
\end{equation*}
Without loss of generality, we assume the strong convexity parameter of $\delta$ is 1, and $\delta(\dot{u}) = 0$. Smooth approximation $\psi_\beta(x)$ with smoothness parameter $\beta > 0$ is defined as 
\begin{equation*}
\psi_\beta(x) = \max_{u \in \mathcal{U}} \left\{ \ip{Bx}{u} - \hat{\phi}(u) - \beta \delta(u) \right\}. 
\end{equation*}
Then, $\psi_\beta$ is well defined, differentiable, convex and smooth. Moreover, it uniformly approximates $\psi$, as it satisfies the following envelop property $\forall x \in \mathcal{X}$:
\begin{equation*}
\psi_\beta(x) \leq \psi(x) \leq \psi_\beta(x) + \beta D_{\mathcal{U}}, 
\end{equation*}
where we denote by $D_{\mathcal{U}} = \max_{u \in \mathcal{U}} \delta(u)$. 
See Theorem~1 in \cite{Nesterov2005} for the proof and  details.

For notational convenience, we restrict ourselves with $g(B\,\cdot\,)$, a Lipschitz continuous function coupled with a linear map, but our findings in this paper directly apply for the general max form. Note that we can write $g(B\,\cdot\,)$ in the max form by choosing $\psi(x) = g(Bx)$, $\hat{\psi}(x) = 0$, and $\hat{\phi}(u)=g^\ast(u)$, Fenchel conjugate of $g$:
\begin{equation*}
g^*(u) = \max_{z} \left\{ \ip{u}{z} - g(z) \right\}.
\end{equation*}
Since $g$ is convex and lower semicontinuous, Fenchel duality holds, and we have $g(Bx) = g^{\ast\ast}(Bx)$. Moreover, Lipschitz continuity assumption on $g$ imposes the boundedness of dual domain. We refer to Lemma~5 by \citet{Dunner2016} for well-known technical details. 

In this work, we specifically focus on the Euclidean prox-functions, $\delta(u) = \frac{1}{2} \norm{u - \dot{u}}^2$. Then, we define $g_\beta$ following the definition of $\psi_{\beta}$ as
\begin{equation*}
g_\beta (Bx) = \max_{u \in \R^q}  \left\{ \ip{Bx}{u} - g^*(u) - \frac{\beta}{2} \norm{u - \dot{u}}^2 \right\}.
\end{equation*}
The argument of this maximization subproblem can be written as $\prox{\beta^{-1}g^\ast}(\beta^{-1}Bx)$, where $\prox{}$ denotes 
\begin{equation*}
\prox{g}(z) = \arg\min_{u} ~ g(u) + \frac{1}{2} \norm{z-u}^2.
\end{equation*}
Hence, following the well-known Moreau decomposition, we can compute the gradient of $g_\beta$ by using
\begin{align*}
\nabla g_\beta (Bx) & = B^\top \prox{\beta^{-1}g^\ast}(\dot{u}  + \beta^{-1}Bx)  \\
& = B^\top \dot{u} + \beta^{-1} B^\top \left(Bx - \prox{\beta g}(\beta\dot{u} + Bx)\right).
\end{align*}

\subsection{Quadratic Penalty} 
Penalty methods often work with unconstrained problems by augmenting the original objective with a penalty function parameterized by a penalty parameter, favoring the constraint. We update this parameter as we progress in the optimization procedure to converge to a solution of the original constrained problem. 

A common and effective proxy is the quadratic penalty, which replaces the affine constraint $Ax \in \mathcal{K}$ by the squared Euclidean distance function, $\frac{\lambda}{2} \dist^2(Ax,\mathcal{K})$, where $\lambda > 0$ is called as the penalty parameter. Surprisingly, quadratic penalty approach is structurally equivalent to a de facto instance of Nesterov smoothing. 

Let us start by writing the Fenchel conjugate of the affine constraint $\iota_{\mathcal{K}}(\,\cdot\,)$,
\begin{equation*}
\iota_{\mathcal{K}}^*(z) = \max_{v \in \mathcal{K}}  \ip{v}{z}.
\end{equation*}
Then, we can write the affine constraint in the max form by choosing $\hat{\psi}(x) = 0$ and $\hat{\phi}(z) = \iota_{\mathcal{K}}^*(z)$, and using the following relation:
\begin{align*}
\iota_{\mathcal{K}}(Ax) 
& = \max_{z} \big\{ \min_{v \in \mathcal{K}} \ip{Ax - v}{z} \big\} \\
& = \max_{z} \left\{ \ip{Ax}{z} - \iota_{\mathcal{K}}^*(z)\right\}.
\end{align*}
Now, by choosing the standard Euclidean prox-function with origin center point, our smooth approximation is
\begin{align*}
{\iota_{\mathcal{K}}}_\beta (Ax) 
& = \max_{z} \left\{ \min_{v \in \mathcal{K}} \ip{Ax - v}{z} - \frac{\beta}{2} \norm{z}^2 \right\} \\
& = \min_{v \in \mathcal{K}} \max_{z} \left\{  \ip{Ax - v}{z} - \frac{\beta}{2} \norm{z}^2 \right\} \\
& =  \frac{1}{2\beta} \dist^2(Ax,\mathcal{K}).
\end{align*}
In summary, we can obtain the quadratic penalty with parameter $\lambda = \beta^{-1}$, by applying Nesterov smoothing procedure to the indicator of an affine constraint. 

Note that the quadratic penalty does not serve as a uniform approximation, because the dual domain is unbounded and the envelop property does not hold. Consequently, the common analysis techniques for  smoothing does not apply for quadratic penalty methods. Nevertheless, one can exploit this structural similarity to design algorithms that universally work on both ends, for composite problems with smoothing friendly non-smooth regularizers, and for problems with affine constraints. In fact, algorithmic design of \citet{Yurtsever2018} implicitly follows this idea. 

Quadratic penalty provides simple and interpretable methods, but with limited practical applicability due to poor empirical performance. To this end, the next subsection reviews augmented Lagrangian methods as an alternative approach. 

\subsection{Augmented Lagrangian}
 
Augmented Lagrangian (AL) methods replace the constraints with a continuous function that promotes feasibility. This function is parametrized by the penalty parameter ($\lambda$, i.e.,\ augmented Lagrangian parameter), and a dual vector $\dot{v}\in\R^p$  (i.e., Lagrange multiplier): 
\begin{align*}
\min_{v\in \mathcal{K}} \left\{ \ip{\dot{v}}{Ax-v} + \frac{\lambda}{2}\norm{Ax-v}^2 \right\}.
\end{align*}

One can motivate augmented Lagrangian penalty from many different point of views. For instance, we can view it as the shifted quadratic penalty, since
\begin{align*}
& \arg\min_{x} f(x) + \min_{v\in \mathcal{K}} \left\{ \ip{\dot{v}}{Ax-v} + \frac{\lambda}{2}\norm{Ax-v}^2 \right\} \\
&\quad = \arg\min_x f(x) + \min_{v\in \mathcal{K}} \, \frac{\lambda}{2}\norm{Ax-v + \tfrac{1}{\lambda}\dot{v}}^2 \\
&\quad = \arg\min_x f(x) +  \frac{\lambda}{2}\dist^2(Ax + \tfrac{1}{\lambda}\dot{v},\mathcal{K}).
\end{align*}

Therefore, we can relate augmented Lagrangian function with Nesterov smoothing in a similar way. To draw this relation, we simply follow the same arguments as in the quadratic penalty case, but this time we use a shifted prox-function $\delta(v) = \frac{1}{2}\norm{v-\dot{v}}^2$:
\begin{align*}
{\iota_{\mathcal{K}}}_\beta (Ax) 
& = \max_{z} \left\{ \min_{v \in \mathcal{K}} \ip{Ax - v}{z} - \frac{\beta}{2} \norm{z - \dot{v}}^2 \right\} \\
& = \min_{v \in \mathcal{K}} \max_{z} \left\{  \ip{Ax - v}{z} - \frac{\beta}{2} \norm{z - \dot{v}}^2 \right\} \\
& =\min_{v\in \mathcal{K}} \left\{ \ip{\dot{v}}{Ax-v} + \frac{1}{2\beta}\norm{Ax-v}^2 \right\}.
\end{align*}

To conclude, augmented Lagrangian formulation is structurally equivalent to an instance of Nesterov smoothing, applied to the indicator of the constraint, with a shifted Euclidean prox-function. The center point of this prox-function corresponds to the dual variable, and penalty parameter corresponds to the inverse smoothness parameter ($\lambda = \beta^{-1}$). 

Once again, this approach does not serve as a uniform approximation, and the common analysis for Nesterov smoothing does not apply for augmented Lagrangian. 

In the next section, we use this basic understanding to design a novel conditional gradient method based on the augmented Lagrangian formulation. 

\section{ALGORITHM}
\label{sec:algorithm}

In this section, we design CGAL for the special case of $g(Bx) = 0$ for the ease of presentation. One can extend CGAL in a straightforward way for the general case, based the discussion in \Cref{sec:prelimineries}, and the analysis techniques in this work and \cite{Yurtsever2018}.

\begin{algorithm}[H]
   \caption{CGAL}
   \label{alg:CGAL}
   {
\begin{algorithmic}
   \State {\bfseries Input:} $x_1 \in \mathcal{X}, ~ y_1 \in \R^p, ~ \lambda_0 > 0$
   \For{$k=1,2, \ldots, $}
	\State $\eta_k = 2/(k +1) \text{~~and~~} \lambda_k = \lambda_0\sqrt{k+1}$.
   	\State $r_k = \proj{\mathcal{K}}\left( Ax_{k} + (1/\lambda_{k})y_k \right) $  
   	\State $v_k =  \nabla f(x_k)  + A^\top y_k + \lambda_{k} A^\top (Ax_k - r_k)$ 
	\State $s_k = \arg\min_{x \in \mathcal{X}}\langle   v_k , x \rangle$
   	\State $x_{k+1} = x_k + \eta_k(s_k - x_k)$  
   	\State $\bar{r}_{k+1} = \proj{\mathcal{K}}\left( Ax_{k+1} + (1/\lambda_{k+1})y_k \right) $  
   	\State $\sigma_{k+1} \gets$ using \eqref{eqn:dual-step-decreasing} or \eqref{eqn:dual-step-constant}
   	\State $y_{k+1} = y_k + \sigma_{k+1} \left( Ax_{k+1} - \bar{r}_{k+1}\right)$ 
   \EndFor
\end{algorithmic}
}
\end{algorithm}

\subsection{Design of CGAL}

Let us introduce the slack variable $r = Ax \in \mathcal{K}$ and define the augmented Lagrangian function as
\begin{align*}\label{eqn:augmented-lagrangian-general}
\mathcal{L}_{\lambda}(x,y) & = f(x) + \min_{r \in \mathcal{K}} \left\{ \ip{y}{Ax - r} + \frac{\lambda}{2}\norm{Ax - r}^2 \right\} \\
& = f(x) - \frac{1}{2\lambda}\norm{y}^2 +  \frac{\lambda}{2}\dist^2\Big(Ax + \frac{1}{\lambda}y,\mathcal{K}\Big).
\end{align*}
where $y \in \R^p$ is the Lagrange multiplier and $\lambda >0$ is the penalty parameter. Clearly $\mathcal{L}_{\lambda}(x,y)$ is a convex $(L_f+\lambda\norm{A}^2)$-smooth function with respect to $x$. 

One CGAL iteration is composed of three basic steps: \\[-2em]
\begin{itemize}[leftmargin=1.5em,itemsep=-0.25em]
\item[$\triangleright$] Primal step (conditional gradient step on $x$),
\item[$\triangleright$] Penalty parameter update (increment $\lambda$),
\item[$\triangleright$] Dual step (proximal gradient step on $y$).
\end{itemize}

\paragraph{Primal step.} 
CGAL is characterized by the conditional gradient step with respect to $\mathcal{L}_{\lambda}(\,\cdot\, ,y)$ on the primal variable. 
At iteration $k$, denoting by
\begin{align*}
r_k = \proj{\mathcal{K}}\Big(Ax_k + \frac{1}{\lambda_{k}}y_k\Big),
\end{align*}
we can evaluate directional gradient $\nabla_x \mathcal{L}_{\lambda_{k}}(x,y_k)$ using
\begin{align*}
\nabla_x \mathcal{L}_{\lambda_k}(x,y_k) = \nabla f(x_k) + A^\top y_k + \lambda_{k} A^\top(Ax_k - r_k).
\end{align*}
Then, we compute linear minimization oracle
\begin{equation*}
s_k = \arg\min_{x \in \mathcal{X}}\ip{ \nabla_x \mathcal{L}_{\lambda_k}(x_k,y_k) }{ x},
\end{equation*}
and we form next iterate ($x_{k+1}$) by combining the current iterate $x_k$ and $s_k$ with CG step-size $\eta_k$. We use the classical step size $\eta_k = 2/(k+1)$ of CG-type methods, but the same guarantees hold for line-search and fully corrective versions.

\paragraph{Penalty parameter update.} 

Penalty methods typically require the penalty parameter to be increased at a certain rate for provable convergence. In contrast, augmented Lagrangian methods can be designed with a fixed penalty parameter, because the saddle point formulation already favors the constraints. Unlike other augmented Lagrangian CG-type methods, we adopt an increasing penalty sequence in CGAL by choosing $\lambda_k = \lambda_0 \sqrt{k+1}$ for some $\lambda_0 > 0$. 

\paragraph{Dual step.} Once $x_{k+1}$ is formed, we update dual variable $y_k$ by a gradient ascent step with respect to $\mathcal{L}_{\lambda}(x, \,\cdot\,)$.
At iteration $k$, we evaluate dual update by
\begin{align*}
y_{k+1} = y_k + \sigma_{k+1} \nabla_y \mathcal{L}_{\lambda_{k+1}}(x_{k+1},y_k).
\end{align*}

To compute $\nabla_y \mathcal{L}_{\lambda_{k+1}}$, we first define
\begin{align*}
\bar{r}_{k+1} = \proj{\mathcal{K}}\big( Ax_{k+1} + \tfrac{1}{\lambda_{k+1}}y_k \big).
\end{align*}
Then, we can use the following formulation:
\begin{align*}
\nabla_y \mathcal{L}_{\lambda_{k+1}}(x_{k+1},y_k) = Ax_{k+1} - \bar{r}_{k+1}.
\end{align*}
Choice of dual step-size is crucial for convergence guarantees. We propose two alternative schemes, with a decreasing or constant bound on the step-size. 

\textit{Decreasing bound on step-size.}
This variant cancels positive quadratic terms in the majorization bounds due to dual updates, with the negative quadratic terms that comes from the penalty parameter updates. Consequently, we choose the largest $\sigma_{k+1} \geq 0$ which satisfies
\begin{equation}\tag{\text{decr.}}
\label{eqn:dual-step-decreasing}
\sigma_{k+1} \leq \tfrac{\lambda_0}{2\sqrt{k+1}} \quad \& \quad \norm{y_{k+1}} \leq D_{\mathcal{Y}_{k+1}} 
\end{equation}
$D_{\mathcal{Y}_{k+1}}$ is a sequence of positive numbers to be chosen, that acts like a dual domain diameter and appears in the final bounds. We will specify a reasonable positive constant $D_{\mathcal{Y}} = D_{\mathcal{Y}_{k+1}}$ in the sequel from the final converges bounds, by matching the factors of the dominating terms. 

\textit{Constant bound on step-size.} We observed significant performance improvements by slightly relaxing the decreasing upper bound on the step-size. To this end, we design this second variant. We do not cancel out additional quadratic terms, but restrict them to be smaller than other dominating terms in the majorization bound. To this end, we choose the largest $\sigma_{k+1} \geq 0 $ which satisfies ($D_{\mathcal{Y}_{k}}$ is similar as in  \eqref{eqn:dual-step-decreasing} case)
\begin{align}
\sigma_{k+1} & \leq \lambda_0 \notag \\
  \norm{y_{k+1}} & \leq D_{\mathcal{Y}_{k+1}}  \tag{\text{const.}} \label{eqn:dual-step-constant} \\
  \quad \sigma_{k+1} \norm{Ax_{k+1} - \bar{r}_{k+1}}^2 & \leq \tfrac{1}{2} \eta_k^2 (L_f + \lambda_{k+1}\norm{A}^2)  D_{\mathcal{X}}^2. \notag
\end{align}
We underline that computation of $\sigma_k$ does not require an iterative line-search procedure, instead it can be computed by simple vector operation both in \eqref{eqn:dual-step-decreasing} and \eqref{eqn:dual-step-constant} variants. As a result, \underline{computational cost of finding $\sigma_k$ is negligible}. 

\subsection{Theoretical Guarantees of CGAL} 

We present convergence guarantees of CGAL in this section, but first we define some basic notions to be used in the sequel and state our main assumptions. 

\paragraph{Solution set.}
We denote a solution of \eqref{eqn:main-template} by $x^\star$, and the set of all solutions by $\mathcal{X}^\star$. Similarly, we denote a solution of the dual problem by $y^\star$, and the set of all solutions by $\mathcal{Y}^\star$. Throughout, we assume that the solution set is nonempty, and that there exists a finite dual solution, \textit{i.e.,} $\min_{y\in\mathcal{Y}^\star} \norm{y} < \infty$. 

\paragraph{$\epsilon$-solution.}
Given an accuracy level $\epsilon > 0$, we call a point $x \in \mathcal{X}$ as an $\epsilon$-solution of \eqref{eqn:main-template} if
\begin{equation*}
f(x) - f^\star \leq \epsilon, \quad \text{and} \quad \dist(Ax,\mathcal{K}) \leq \epsilon.
\end{equation*}
We call $f(x) \!- \! f^\star$ as the objective residual and $\dist(Ax,\mathcal{K})$ as the feasibility gap. Note that the convergence of objective residual alone is not enough to approximate the solution, since the iterates are non-feasible and $f(x)\! -\! f^\star$ can take negative values. 

\paragraph{Strong duality.} 
We assume Slater's condition, which is a  sufficient condition for strong duality: 
\begin{equation*}
\mathrm{relint}(\mathcal{X} \times \mathcal{K}) ~\cap~ \{(x,r) \in \mathrm{dom}(f) \times \R^d : Ax = r \} \neq \emptyset,
\end{equation*}
where $\mathrm{relint}$ means relative interior. Strong duality is a common assumption for primal-dual methods.

\begin{theorem}\label{thm:cgal-decreasing}
Sequence $x_k$ generated by CGAL with dual step-size conditions \eqref{eqn:dual-step-constant} satisfies:
\begin{align*}
f(x_{k}) - f^\star & \geq -\norm{y^\star} ~ \mathrm{dist}(Ax_k,\mathcal{K}) \\
f(x_{k}) - f^\star & \leq 4 D_{\mathcal{X}}^2 \left( \frac{L_f}{k}   + \frac{\lambda_0 \norm{A}^2}{\sqrt{k}} \right) + \frac{D_{\mathcal{Y}_{k}}^2}{2\lambda_0 \sqrt{k}} \\
\mathrm{dist}(Ax_k,\mathcal{K})  & \leq \frac{2/\lambda_0}{\sqrt{k}} \! \left( \! \frac{D_{\mathcal{Y}_k}}{2} \!+\! \norm{y_{k} - y^\star} \! + \! \sqrt{2C_0\lambda_0 D_{\mathcal{X}}^2 } \right)
\end{align*}
where $C_0 = L_f + \|A\|^2 \lambda_{0}$. We can also bound $\norm{y_{k} - y^\star}$ using triangle inequality. Considering the bounds, it is reasonable to choose $D_{\mathcal{Y}}$ proportional to $D_{\mathcal{X}}\norm{A}\lambda_0$.

Sequence $x_k$ generated by CGAL with dual step-size conditions \eqref{eqn:dual-step-decreasing} satisfies similar guarantees as \eqref{eqn:dual-step-constant}, with the factor of $1/2$ for all terms involving $D_{\mathcal{X}}^2$.
\end{theorem}

We omit design variants of CGAL with line-search and fully corrective updates, covered by our theory.  

\subsection{Extension for Composite Problems}\label{sec:composite-extension}

One can extend CGAL in a straightforward way for composite problems based on the discussions in \Cref{sec:prelimineries}. For this, we simply need to define the sum of two non-smooth terms: $G(Ax,Bx) = \iota_{\mathcal{K}}(Ax) + g(Bx).$ Then, we consider smooth approximation of this term with smoothness parameter $\beta = 1/\lambda$ and prox-function $\delta(u,v) = \frac{1}{2}\norm{u-\dot{u}}^2 + \frac{1}{2}\norm{v-\dot{v}}^2$:
\begin{align*}
G_{\beta}(Ax,Bx) = {\iota_{\mathcal{K}}}_{\beta}(Ax) + g_{\beta}(Bx)
\end{align*}
Gradient of $G_{\beta}$ can be written as the sum of individual gradient terms. 
Then, CGAL applies simply by adding one more dual variable, $z \in \R^q$ and changing $v_k$ as 
\begin{align*}
v_k ~\gets &  ~\nabla  f(x_k)   + A^\top y_k + \lambda_{k} A^\top (Ax_k - r_k) \\ 
& + B^\top z_k + \lambda_{k} B^\top (Bx - \prox{\frac{1}{\lambda_k} g}(Bx + \tfrac{1}{\lambda_k} z_k)).
\end{align*}
One can keep $z$ fixed as in \cite{Yurtsever2018}, or update it similar to $y$ using $\prox{\lambda g^*}$. Both cases guarantees $\mathcal{O}(1/\sqrt{k})$ rates in feasibility gap $\dist(Ax,{\mathcal{K}})$ and in objective residual $f(x) + g(Bx) - f(x^\star) - g(Bx^\star)$.

\section{RELATED WORK}
\label{sec:related-work}

The majority of convex methods for solving \eqref{eqn:main-template} are based on computationally challenging oracles, which can be some second order oracle as in interior point methods, a projection step (onto $\mathcal{X}$) as in operator splitting methods, or a constrained proximal-oracle as in classical primal-dual methods. We refer to \cite{Wright1997},\cite{Komodakis2015},\cite{Ryu2016} and the references therein for these classical approaches. In the rest of this section we focus on the optimization methods which applies \eqref{eqn:main-template} or some of its subclasses by leveraging the linear minimization oracle. 

{\bf \citet{Lan2014}} introduces a conditional gradient method for non-smooth minimization over a convex compact domain, based on Nesterov smoothing. To the best of our knowledge, this is the first attempt to combine Nesterov smoothing and conditional gradient approach. This method does not apply for our general problem template, in particular to problems with affine constraints, since it relies on the boundedness of the dual domain and the uniform approximation property. 

{\bf \citet{Yurtsever2015}} present the universal primal-dual method (UPD), a primal-dual subgradient approach for solving convex minimization problems under affine constraints. Main template of UPD is fairly different than \eqref{eqn:main-template}, it does not have the non-smooth term $g(Bx)$ instead it assumes H\"older smoothness in the dual space. The method does not directly work with $\mathrm{lmo}$'s, but it leverages the so-called sharp operators with comparable computational complexity to $\mathrm{lmo}$'s under some specific problem settings. In particular, for standard SDP's with linear cost function, sharp operator becomes the same as $\mathrm{lmo}$

UPD adopts the inexact line-search strategy introduced by \citet{Nesterov2015}. This strategy requires the input of target accuracy $\epsilon$, and UPD is guaranteed to converge only up to $\epsilon$ accuracy, i.e., it guarantees $f(x) - f^\star \leq \mathcal{O}(1/\sqrt{k}) + \epsilon$. Practical performance of this method heavily depends on this parameter: Choosing $\epsilon$ small causes step-sizes to be too small. The best value of $\epsilon$ is typically around $1/10$th and $1/100$th of the optimal value $|f^\star|$, and this method is difficult to tune unless optimal value is roughly known. 

{\bf \citet{Lan2016}} propose the conditional gradient sliding method (CGS). This method is based on an inexact version of accelerated gradient method by \citet{Nesterov1983}, where the projection oracle is approximated by CGM. CGS is originally proposed for smooth minimization over a convex and compact domain, but the results are generalized for smoothing friendly non-smooth functions in Section~4, following the same approach as \citet{Lan2014}. Note that this generalization directly follows the standard approach of Nesterov smoothing, and it does not apply for affine constraints. 

{\bf \citet{yen2016convex}} proposes the greedy direction method of multipliers (GDMM), a CGM variant for minimizing a linear objective function over an intersection of polytopes. GDMM relies on a consensus reformulation over cartesian product of these polytopes, and the consistency constraint is incorporated by the augmented Lagrangian. This method is further explored in structural support vector machine \cite{Yen2016svm} and maximum-a-posteriori inference \cite{Huang2017} problems. Nevertheless, as raised later on by \cite{Gidel2018}, there are technical issues in the analysis which do not admit a trivial fix. We refer to Section~B.1 in \cite{Gidel2018} for more details. 

{\bf \citet{Gidel2018}} propose an augmented Lagrangian framework for convex splitting problem (FW-AL). Similar to CGAL, this method is characterized by one CGM step on $\mathcal{L}_{\lambda}(\cdot,y_k)$ followed by one dual gradient ascent step on $\mathcal{L}_{\lambda}(x_{k+1},\cdot)$, but their penalty parameter $\lambda$ is fixed. Their method is specific for $Ax=0$ (i.e., splitting), but it can be applied to $Ax = b$ case using a product space technique. The analysis of FW-AL relies on the error bounds (see Theorem~1 in \cite{Gidel2018} for the conditions, and \cite{Bolte2017} for more details about error bounds). Their step-size $\sigma_{k+1}$ depends on the error bound constant $\alpha$ as $\sigma_{k+1} = \frac{2\sigma_0}{k+2}$ with $\sigma_0 \leq \min\{\frac{2}{\lambda}, \frac{\alpha^2}{2\delta} \}$. Hence, $\sigma_0$ is a tuning parameter, and the method has guarantees only if it is chosen small enough. Note that $\alpha$ is not only unknown, it can be also arbitrarily small. 

{\bf \citet{Liu2018}} introduce an inexact augmented Lagrangian method (IAL), where the Lagrangian subproblems are approximated by CGM up to a prescribed accuracy level, say $\varepsilon_k = \varepsilon_0/k$ for some $\varepsilon_0 > 0$ to be tuned. This results in a double-loop algorithm, where each iteration consists multiple CGM iterations until the following condition is satisfied:
\begin{align*} 
\max_{x \in \mathcal{X}} \ip{\nabla f(x_{k+1}) + A^\top y_k + \lambda A^\top (Ax_{k+1} - b)}{x} \leq \varepsilon_k.
\end{align*}
Then, the algorithm takes a dual gradient ascent step. 

IAL provably generates an $\epsilon$-solution after $\mathcal{O}(1/\epsilon^2)$ outer iterations, by choosing the penalty parameter $\lambda$ appropriately (proportional to $1/\sqrt{\epsilon}$). This method, however, requires multiple $\mathrm{lmo}$ calls at each iteration. Since the number of $\mathrm{lmo}$ calls is bounded by $\lceil{6 L_f D_{\mathcal{X}}^2 / \varepsilon_k} \rceil - 2$ (see Theorem~{2.2} in \cite{Liu2018}), this results in $\mathcal{O}(1/\epsilon^4)$ calls of $\mathrm{lmo}$. Note that this is much worse than $\mathcal{O}(1/\epsilon^2)$ calls required by our method. 

{\bf \citet{Yurtsever2018}} present a conditional gradient type method (HCGM) for \eqref{eqn:main-template}. This method relies on the quadratic penalty approach to handle affine constraints. HCGM guarantees $\mathcal{O}(1/\sqrt{k})$ convergence rate both in the objective and the feasibility gap similar to CGAL. Note however, as explained in \Cref{sec:prelimineries}, penalty methods typically performs with the worst case rates. We can indeed observe this in numerical experiments of \citet{Yurtsever2018}, that the empirical convergence rate is also $\mathcal{O}(1/\sqrt{k})$. We demonstrated this also in our experiments in \Cref{sec:experiments}.

\section{NUMERICAL EXPERIMENTS}
\label{sec:experiments}
\begin{figure*}[ht]
\centering
\includegraphics[width=0.49\textwidth]{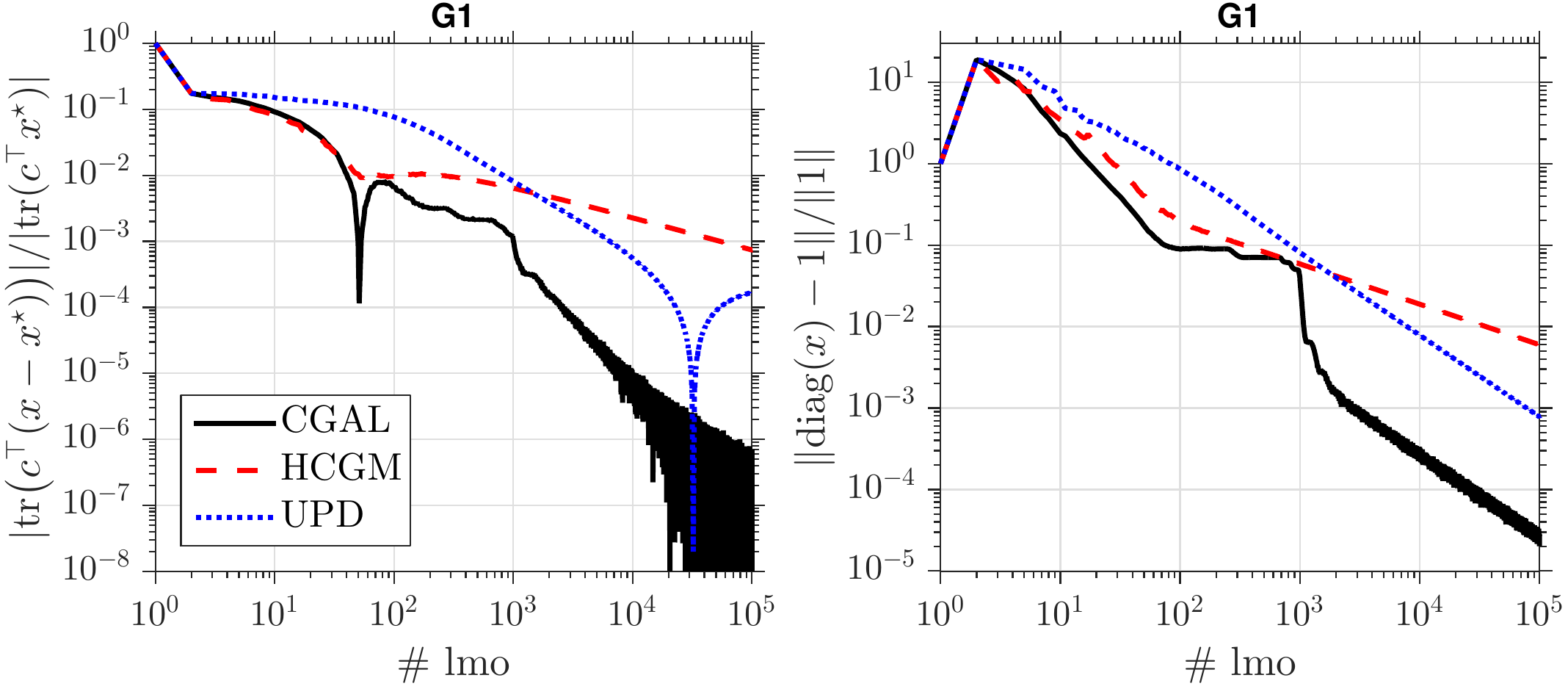}
\includegraphics[width=0.49\textwidth]{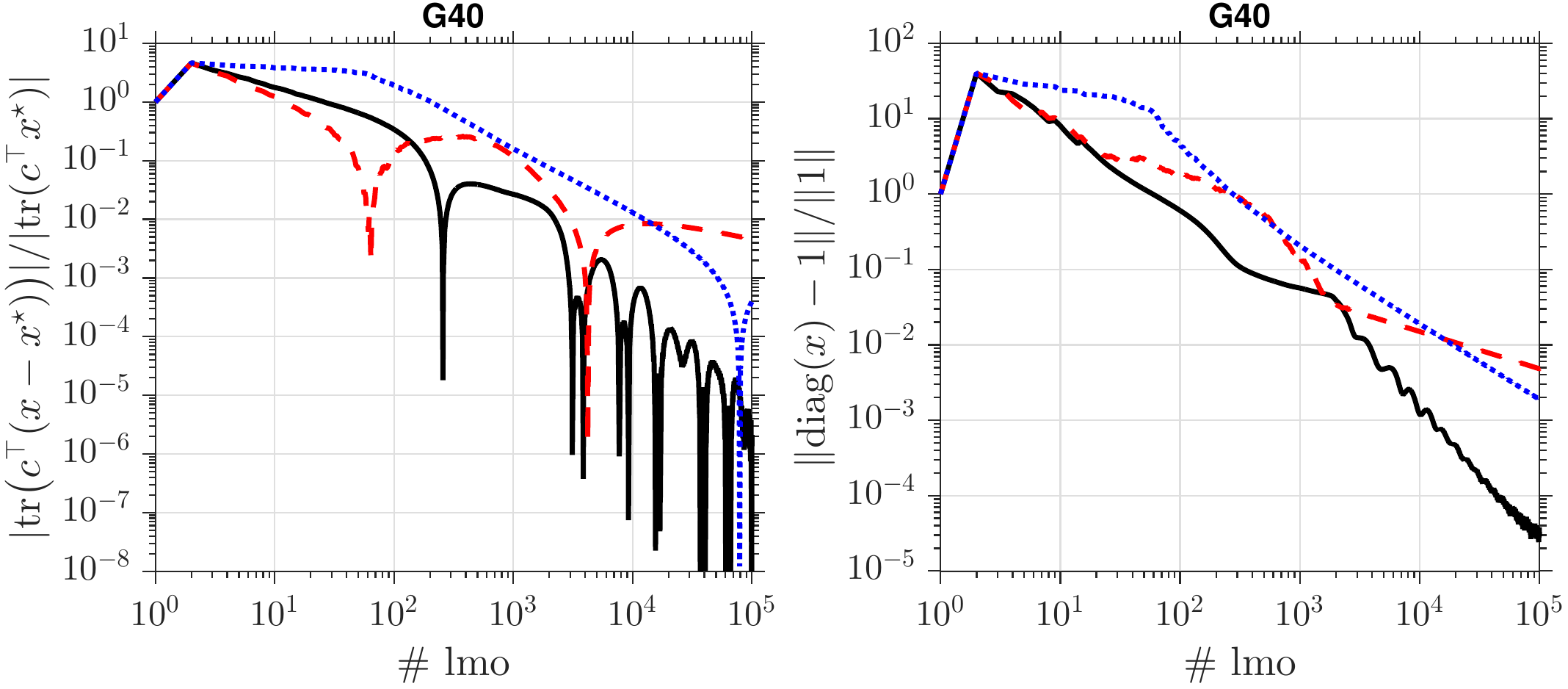}\\[-1em]
        \caption{\small{Empirical comparison of CGAL, HCGM and UPD with max-cut problem setup.}}
        \label{fig:maxcut-gsets}
        \vspace{-0.5em}
\end{figure*}

This section presents the numerical evidence to demonstrate empirical superiority of CGAL, based on max-cut, clustering and generalized eigenvector problems. 

We compared CGAL against UPD and HCGM from \Cref{sec:related-work}. This choice is based on the practicality of the algorithms: FW-AL and IAL have $2$ tuning parameters each, and it is very difficult to use these methods in medium or large scale experiments. CGAL and HCGM has initial penalty parameter $\lambda_0$, and UPD has accuracy parameter $\epsilon$ to be tuned. We tuned all these parameters by bisection with factor $10$, until the method with the chosen parameter outperforms itself both with $10$th and $1/10$th of the parameter. Although CGAL with \eqref{eqn:dual-step-decreasing} performed better than HCGM in all instances we tried, CGAL with \eqref{eqn:dual-step-constant} uniformly outperformed \eqref{eqn:dual-step-decreasing} and HCGM. Hence in this section we focus on CGAL with \eqref{eqn:dual-step-constant}. 

Note that the computational cost of all algorithms are dominated by $\mathrm{lmo}$, hence we provide plots with number of $\mathrm{lmo}$ calls on the x-axis which is roughly proportional to computation time. 

\subsection{Max-cut}\label{sec:experiment-max-cut}

\begin{figure}[ht!]
\centering
\includegraphics[width=\linewidth]{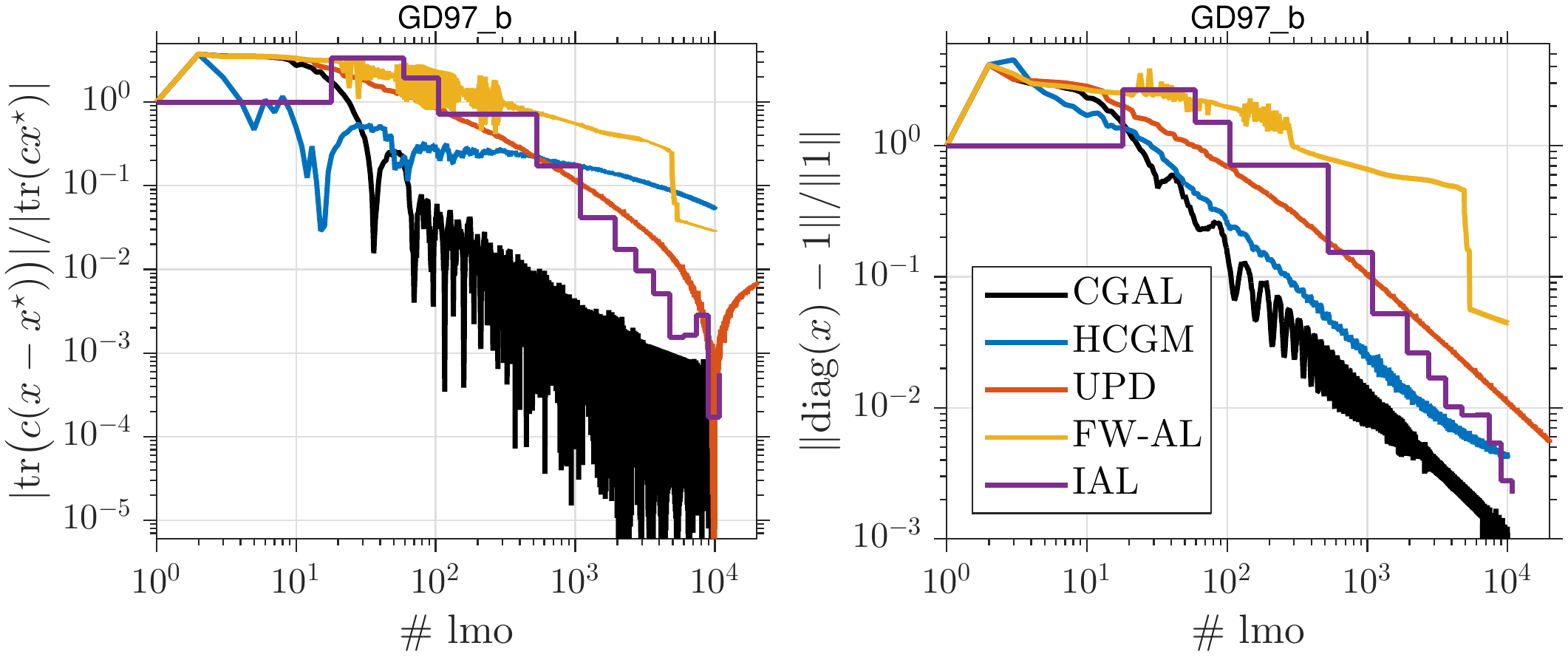}\\[-0.5em]
        \caption{\small{Overview of empirical performance of various methods in \Cref{sec:related-work} for max-cut problem.}}
        \label{fig:maxcut-small}
\end{figure}

Maximum cut is an NP-Hard combinatorial problem from computer science. Denoting the symmetric $n \times n$ graph Laplacian matrix of a graph by $c$, this problem can be relaxed as \cite{Goemans1995}:
\begin{equation*}
\begin{aligned}
&\underset{x}{\text{maximize}} 
& & \tfrac{1}{4} \trace(c x) \\
&\text{subject to} 
& & \trace (x) = n, ~~ x \in \mathbb{S}_+^n ~~  \& ~~ \text{diag}(x) = 1.
\end{aligned}
\end{equation*}
Tuning all methods from \Cref{sec:related-work} require substantial computational effort, especially since some of these methods have multiple tuning parameters. To this end, we first consider a small scale max-cut instance where we compare against all of these methods. In this setup we use GD97\_b dataset \cite{GD97b}, which corresponds to a $47 \times 47$ graph. In \Cref{fig:maxcut-small}, we present the performance of each method with the best parameter choice obtained after extensive search. We also provide the performance with all trials of each algorithm in the supplements, also with some other variants of the methods. 

Next, we move from the toy case to medium scale examples, and we compare CGAL with UPD and HCGM for max-cut with G1 ($800 \times 800$) and  G40 ($2000 \times 2000$) datasets \cite{GSet}. We compile the results of these tests in \Cref{fig:maxcut-gsets}, from which we observe that HCGM converges with $\mathcal{O}(1/\sqrt{k})$ (which is the worst case bound) while CGAL achieves faster than $\mathcal{O}(1/k)$ rate. Note that the sudden drop of UPD on the objective residual plots towards the end is not an increase of rate, it is simply the sign flip of $f(x) - f^\star$ which typically happens just before the saturation of UPD. 

\subsection{k-means Clustering}

We consider a test setup with SDP formulation of model-free k-means clustering by \cite{Peng2007}:
\begin{equation*}
\begin{aligned}
&\underset{x}{\text{minimize}} 
& & \trace(c x) \\
&\text{subject to} 
& & x1 = 1, ~~ x\geq 0, ~~ x \in \mathbb{S}_+^n ~~  \& ~~ \trace(x) = \alpha.
\end{aligned}
\end{equation*}
where $\alpha$ is the number of clusters, and $c$ is $n \times n$ Euclidean distance matrix. We denote by $1$ the vector of ones, hence $x1 = 1$ and $x \geq 0$ together enforce each row to be on the unit simplex. Same applies for columns due to symmetry. We cast this problem into \eqref{eqn:main-template} by choosing $f(x) = \trace(c x)$, $\mathcal{X} = \{x: x \in \mathbb{S}_+^n , \trace(x) = \alpha \}$, $A:\mathbb{S}_+^n \to \R^n \times \R^{n\times n}$ maps $x \to (x1 , x)$, and finally $\mathcal{K} = \{1\} \times \R_+^{n\times n}$.

\begin{figure}[ht!]
\centering
\includegraphics[width=\linewidth]{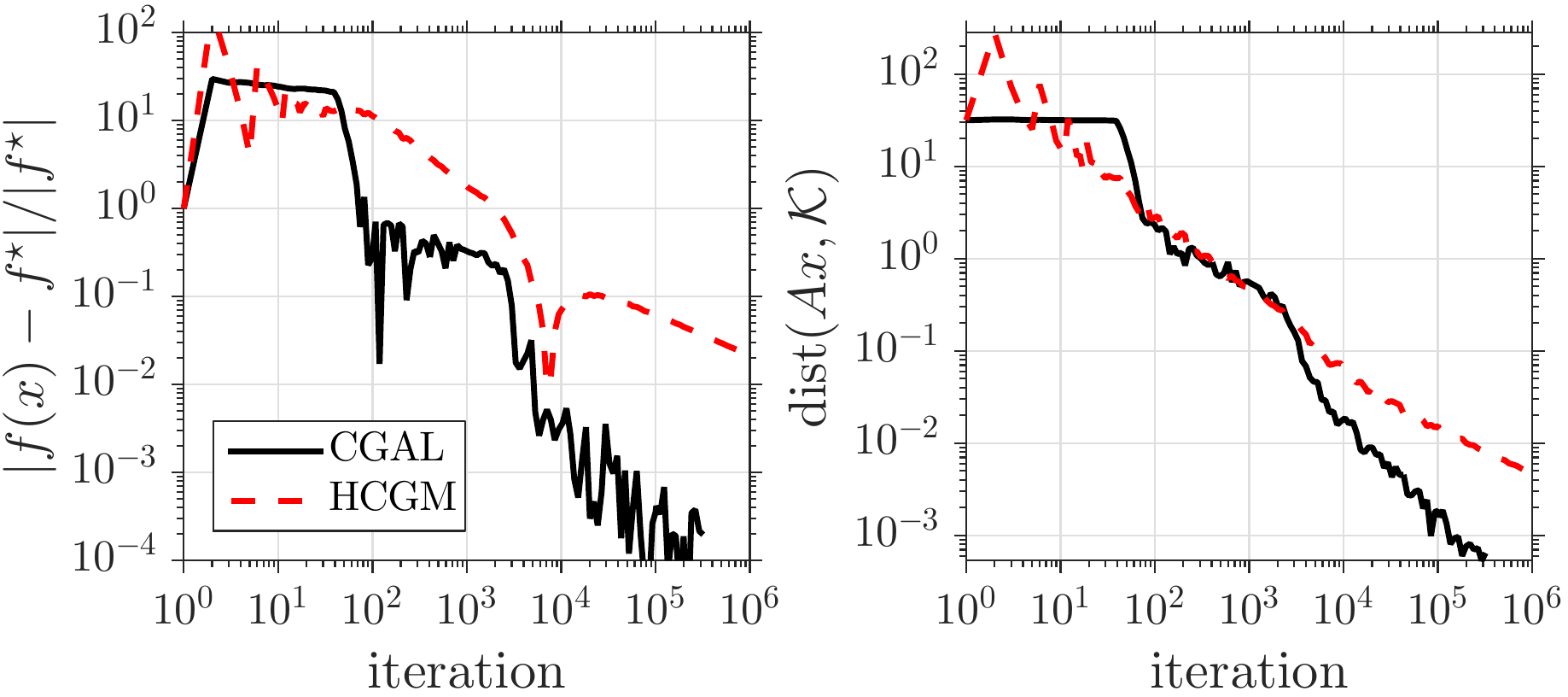}\\[-0.5em]
        \caption{\small{Objective residual and feasibility gap for k-means clustering with preprocessed MNIST dataset.}}
        \label{fig:clustering}
\end{figure}

\begin{figure*}[ht!]
\centering
\includegraphics[width=0.49\textwidth]{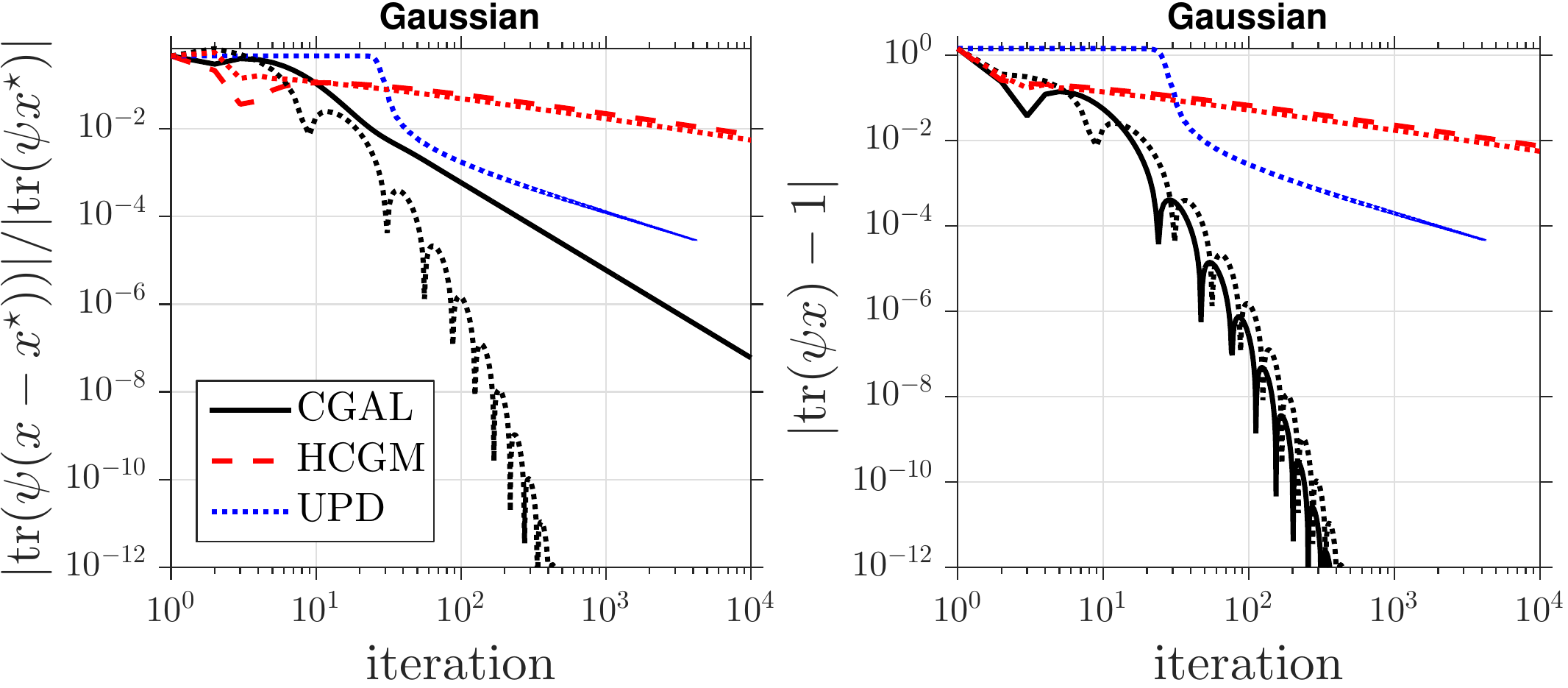}
\includegraphics[width=0.49\textwidth]{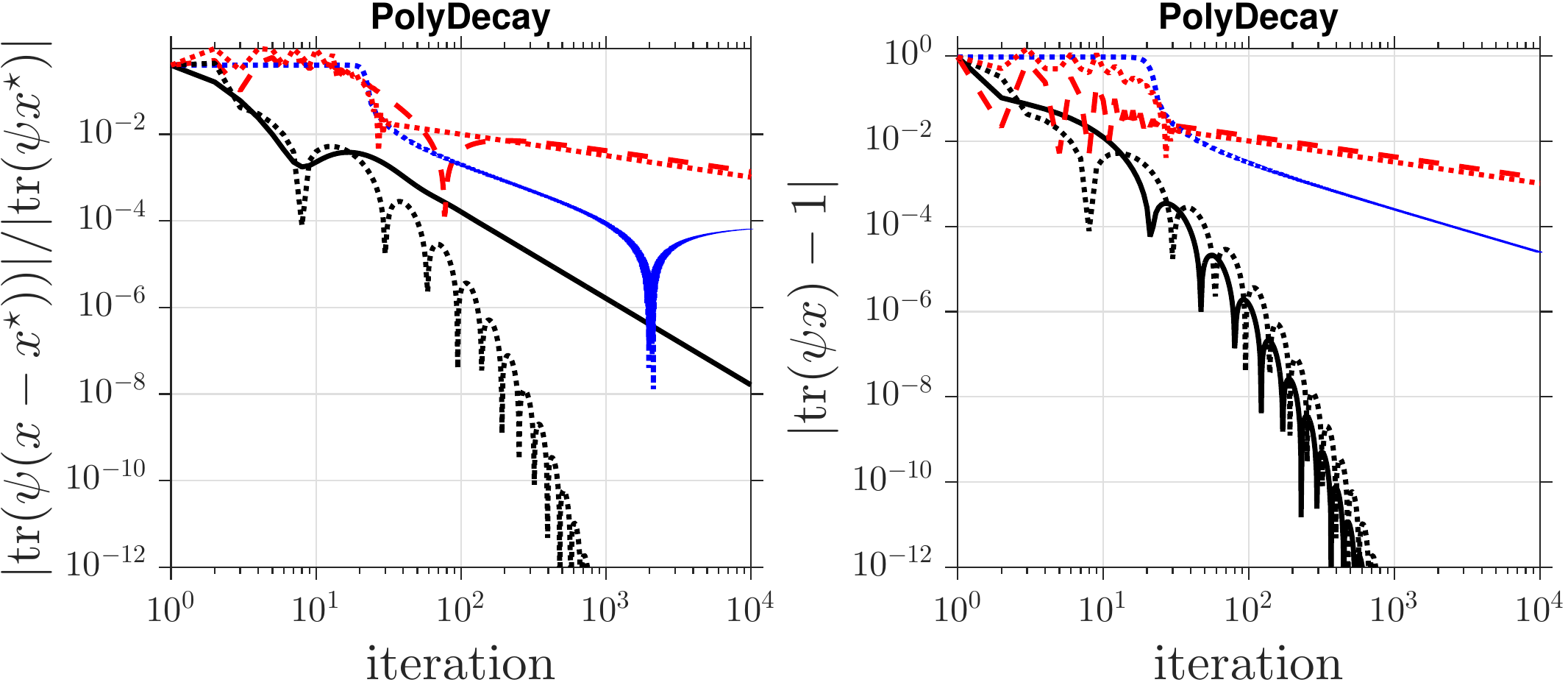} \\[0.5em]
\includegraphics[width=0.49\textwidth]{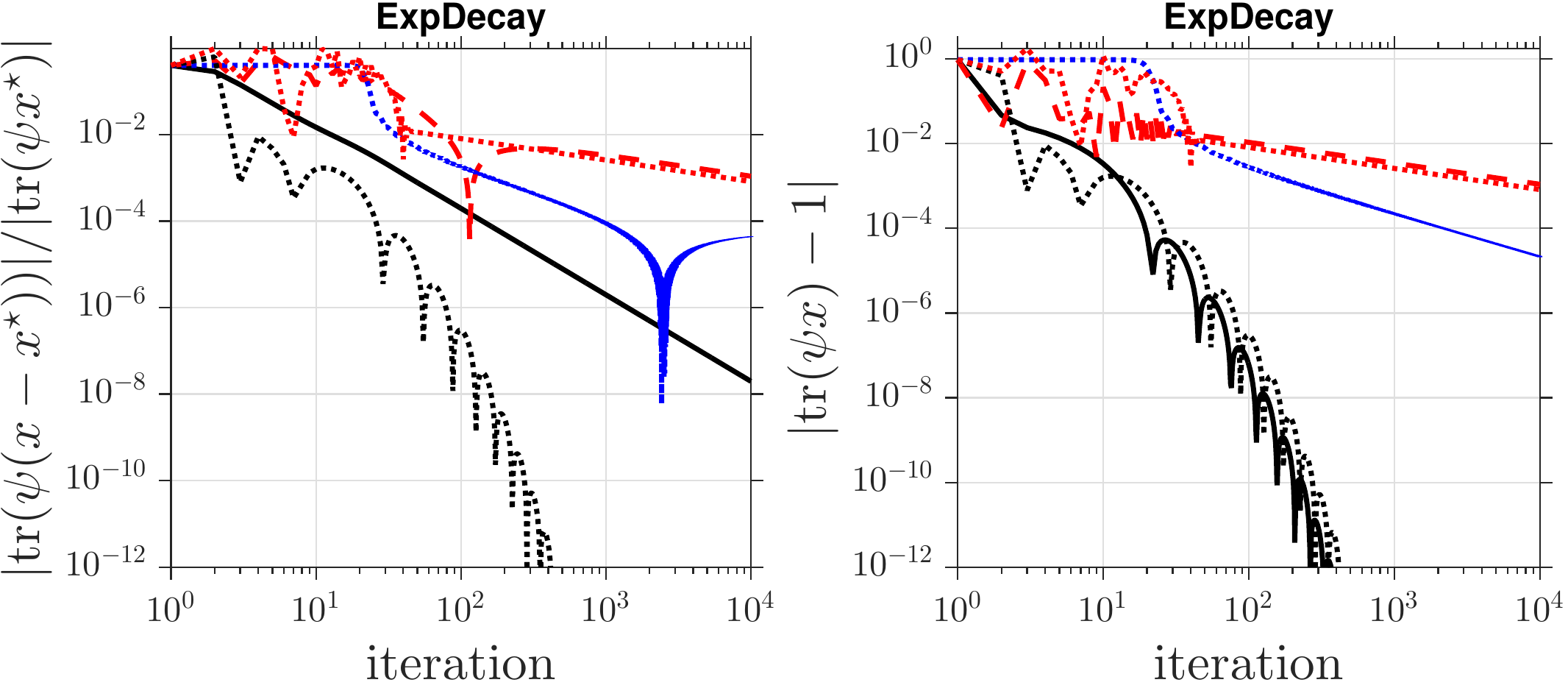}
\includegraphics[width=0.49\textwidth]{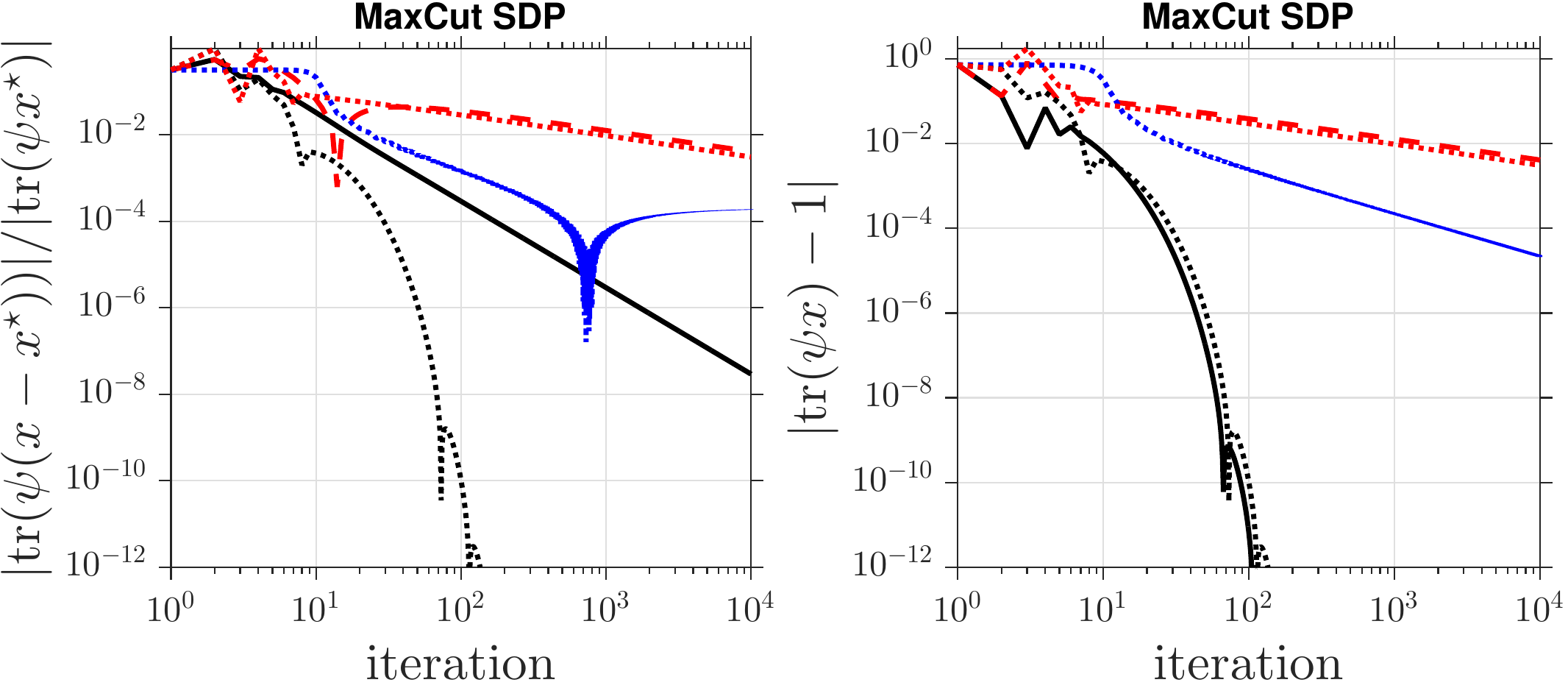}\\[-1em]
        \caption{\small{Empirical comparison of CGAL \& HCGM for solving generalized eigenvector problem with 4 different synthetic setups. Dotted lines present objective residual and feasibility gap of the atoms chosen by linear minimization oracle ($s_k$).}}
        \label{fig:geneig}
\end{figure*}

We use the same setup as in \cite{Yurtsever2018}, which is designed and published online by \citet{Mixon2017}. This setup contains a $1000 \times 1000$ dimensional dataset generated by sampling and preprocessing the MNIST dataset \cite{MNIST} using a one-layer neural network. Further details on this setup and the dataset can be found in \cite{Mixon2017}. 

In \Cref{fig:clustering}, we observe once again that CGAL outperforms HCGM, achieving $\mathcal{O}(1/k)$ empirical convergence rate. In this problem instance, we failed to tune UPD, even with the knowledge of $f^\star$. After extensive analysis and tests, we concluded that UPD has an implicit tuning parameter. It is possible to choose different accuracy terms for objective and feasibility in UPD, as also noted by the authors, simply by scaling the objective function with a constant. Performance of UPD heavily depends on this scaling in addition to tuning accuracy parameter, hence we omit UPD. 

\subsection{Generalized Eigenvector Problem}

We consider SDP relaxation of the generalized eigenvector problem from \citet{Boumal2018}:
\begin{equation*}
\begin{aligned}
&\underset{x}{\text{maximize}} 
& & \trace(\phi x) \\
&\text{subject to} 
& & \trace(x) \leq \alpha, ~~ X \in \mathbb{S}_+^n ~~  \& ~~ \trace(\psi x) = 1
\end{aligned}
\end{equation*}
where $\phi$ and $\psi$ are symmetric matrices of size $n \times n$, and $\alpha >0$ is a model parameter. In this problem, we consider synthetic setups, where we generate $\psi$ iid Gaussian, and consider $4$ different cases for $\phi$:\vspace{-0.75em}
\begin{itemize}[leftmargin=1em,itemsep=-0.15em]
\renewcommand{\labelitemi}{$\circ$}
\item \texttt{Gaussian} - $\phi$ generated by taking symmetric part of $10^3 \times 10^3$ iid Gaussian matrix
\item \texttt{PolyDecay} - $\phi$ generated by randomly rotating  $\mathrm{diag}(1^{-i},2^{-i},\ldots,1000^{-i})$ ($i=1$)
\item \texttt{ExpDecay} - $\phi$ generated by randomly rotating $\mathrm{diag}(10^{-i},10^{-2i},\ldots,10^{-1000i})$ ($i = 0.025$)
\item \texttt{MaxCut SDP} - $\phi$ is a solution of a maxcut SDP with G40 dataset ($2000 \times 2000$) \vspace{-0.75em}
\end{itemize}

This problem highlights an important observation under various data models, which empirically explains why CGAL outperforms the base method HCGM. Note that this SDP problem provably has a rank-$1$ solution, and if $\alpha$ is tuned to its exact value, the solution is an extreme point of the domain. In this scenario, we might expect $\mathrm{lmo}$ to pick the solution itself, or other close extreme points, if the problem formulation is well-conditioned. Recall that CGAL updates the dual variable, which corresponds to the center point of a quadratic penalty, with the expectation of better adaptation to the problem geometry. In \Cref{fig:geneig}, we provide an empirical evidence of this adaptation, where dotted lines corresponds to extreme points chosen by $\mathrm{lmo}$. Unsurprisingly, $s_k$ converges quickly (with linear rates) under different scenarios for CGAL, while we do not observe the same behavior in HCGM or UPD (we omit lmo outputs of UPD in figure which do not converge).

\section{CONCLUSIONS}
\label{sec:conclusions}

CGAL can be viewed as a natural extension of HCGM of \cite{Yurtsever2018} going from quadratic penalty to an augmented Lagrangian formulation. CGAL retains the strong theoretical guarantees of HCGM as well as (nearly) the same per-iteration complexity while exhibiting significantly superior empirical performance ($\mathcal{O}(1/k)$ in practice vs $\mathcal{O}(1/\sqrt{k})$ in theory). In stark contrast to the existing methods that apply to the SDP templates, CGAL does not require strong assumptions on the problem geometry for the optimal convergence rates, and it comes from a simple analysis with interpretable bounds from which we can derive analytical dual step-size rules. Considering recent developments on the storage optimal convex optimization methods based on conditional gradients \cite{Yurtsever2016}, CGAL might be the key step for designing fast convex optimization methods solving for huge scale SDP's. Further exploration of CGAL for this specific problem setup is left for future research in addition to proving the faster convergence rate, which requires stronger analytical assumptions that should also hold for SDPs. 

\subsubsection*{Acknowledgements}

$^\dagger$This work was supported by the Swiss National Science Foundation (SNSF) under grant number $200021\_178865 / 1$. This project has received funding from the European Research Council (ERC) under the European Union's Horizon $2020$ research and innovation programme (grant agreement no $725594$ - time-data). $^\ddagger$This work was supported by a public grant as part of the Investissement d'avenir project, reference ANR-$11$-LABX-$0056$-LMH, LabEx LMH, in a joint call with PGMO. 

\bibliography{../../common/bibliography.bib}

\begin{thebibliography}{31}
\providecommand{\natexlab}[1]{#1}
\providecommand{\url}[1]{\texttt{#1}}
\expandafter\ifx\csname urlstyle\endcsname\relax
  \providecommand{\doi}[1]{doi: #1}\else
  \providecommand{\doi}{doi: \begingroup \urlstyle{rm}\Url}\fi

\bibitem[Batagelj and Mrvar()]{GD97b}
V.~Batagelj and A.~Mrvar.
\newblock Pajek datasets, http://vlado.fmf.uni-lj.si/pub/networks/data/.

\bibitem[Bolte et~al.(2017)Bolte, Nguyen, Peypouquet, and Suter]{Bolte2017}
J.~Bolte, T.~P. Nguyen, J.~Peypouquet, and B.~W. Suter.
\newblock From error bounds to the complexity of first-order descent methods
  for convex functions.
\newblock \emph{Math. Program.}, 165\penalty0 (2):\penalty0 471--507, 2017.

\bibitem[Boumal et~al.(2018)Boumal, Voroninski, and Bandeira]{Boumal2018}
N.~Boumal, V.~Voroninski, and A.~Bandeira.
\newblock Deterministic guarantees for {B}urer--{M}onteiro factorizations of
  smooth semidefinite programs.
\newblock arXiv:1804.02008v1, 2018.

\bibitem[Cand\`es et~al.(2013)Cand\`es, Strohmer, and Voroninski]{Candes2011}
E.~J. Cand\`es, T.~Strohmer, and V.~Voroninski.
\newblock {P}hase{L}ift: Exact and stable signal recovery from magnitude
  measurements via convex programming.
\newblock \emph{Communications on Pure and Applied Math.}, 66\penalty0
  (8):\penalty0 1241--1274, 2013.

\bibitem[D\"unner et~al.(2016)D\"unner, Forte, Tak\'ac, and Jaggi]{Dunner2016}
C.~D\"unner, S.~Forte, M.~Tak\'ac, and M.~Jaggi.
\newblock Primal--dual rates and certificates.
\newblock In \emph{Proc. $33$rd Int. Conf. Machine Learning}, 2016.

\bibitem[Frank and Wolfe(1956)]{FrankWolfe1956}
M.~Frank and P.~Wolfe.
\newblock An algorithm for quadratic programming.
\newblock \emph{Naval Research Logistics Quarterly}, 3:\penalty0 95--110, 1956.

\bibitem[Gidel et~al.(2018)Gidel, Pedregosa, and Lacoste-Julien]{Gidel2018}
G.~Gidel, F.~Pedregosa, and S.~Lacoste-Julien.
\newblock {F}rank-{Wolfe} splitting via augmented {L}agrangian method.
\newblock In \emph{Proc. $21$st Int. Conf. Artificial Intelligence and
  Statistics (AISTATS)}, 2018.

\bibitem[Goemans and Williamson(1995)]{Goemans1995}
M.~X. Goemans and D.~P. Williamson.
\newblock Improved approximation algorithms for maximum cut and satisfiability
  problems using semidefinite programming.
\newblock \emph{Journal of the ACM}, 43\penalty0 (6):\penalty0 1115--1145,
  1995.

\bibitem[Hazan and Kale(2012)]{Hazan2012}
E.~Hazan and S.~Kale.
\newblock Projection--free online learning.
\newblock In \emph{Proceedings of the $29$th International Conference on
  Machine Learning}, 2012.

\bibitem[Huang et~al.(2017)Huang, Yen, Zhang, Huang, Ravikumar, and
  Dhillon]{Huang2017}
X.~Huang, I.~E.-H. Yen, R.~Zhang, Q.~Huang, P.~Ravikumar, and I.~S. Dhillon.
\newblock Greedy direction method of multiplier for {MAP} inference of large
  output domain.
\newblock In \emph{Proc. $20$th Int. Conf. Artificial Intelligence and
  Statistics (AISTATS)}, 2017.

\bibitem[Jaggi(2013)]{Jaggi2013}
M.~Jaggi.
\newblock Revisiting {F}rank--{W}olfe: {P}rojection--free sparse convex
  optimization.
\newblock In \emph{Proc. $30$th Int. Conf. Machine Learning}, 2013.

\bibitem[Komodakis and Pesquet(2015)]{Komodakis2015}
N.~Komodakis and J.-C. Pesquet.
\newblock Playing with duality: An overview of recent primal-dual approaches
  for solving large-scale optimization problems.
\newblock \emph{IEEE Signal Process. Mag.}, 32\penalty0 (6):\penalty0 31--54,
  2015.

\bibitem[Lan(2014)]{Lan2014}
G.~Lan.
\newblock The complexity of large--scale convex programming under a linear
  optimization oracle.
\newblock arXiv:1309.5550v2, 2014.

\bibitem[Lan and Zhou(2016)]{Lan2016}
G.~Lan and Y.~Zhou.
\newblock Conditional gradient sliding for convex optimization.
\newblock \emph{SIAM J. Optim.}, 26\penalty0 (2):\penalty0 1379--1409, 2016.

\bibitem[LeCun and Cortes()]{MNIST}
Y.~LeCun and C.~Cortes.
\newblock {MNIST} handwritten digit database, Accessed: Jan. 2016 .
\newblock URL \url{http://yann.lecun.com/exdb/mnist/}.

\bibitem[Liu et~al.(2018)Liu, Liu, and Ma]{Liu2018}
Y.-F. Liu, X.~Liu, and S.~Ma.
\newblock On the non-ergodic convergence rate of an inexact augmented
  {L}agrangian framework for composite convex programming.
\newblock \emph{to appear in Mathematics of Operations Research}, 2018.

\bibitem[Mixon et~al.(2017)Mixon, Villar, and Ward]{Mixon2017}
D.~G. Mixon, S.~Villar, and R.~Ward.
\newblock Clustering subgaussian mixtures by semidefinite programming.
\newblock \emph{Information and Inference: A Journal of the IMA}, 6\penalty0
  (4):\penalty0 389--415, 2017.

\bibitem[Nesterov(1987)]{Nesterov1983}
Y.~Nesterov.
\newblock A method of solving a convex programming problem with convergence
  rate {$O( 1 / k^2 )$}.
\newblock \emph{Soviet Mathematics Doklady}, 27\penalty0 (2):\penalty0
  372--376, 1987.

\bibitem[Nesterov(2005{\natexlab{a}})]{Nesterov2005}
Y.~Nesterov.
\newblock Smooth minimization of non-smooth functions.
\newblock \emph{Math. Program.}, 103:\penalty0 127--152, 2005{\natexlab{a}}.

\bibitem[Nesterov(2005{\natexlab{b}})]{Nesterov2005eg}
Y.~Nesterov.
\newblock Excessive gap technique in nonsmooth convex minimization.
\newblock \emph{SIAM J. Optim.}, 16\penalty0 (1):\penalty0 235--249,
  2005{\natexlab{b}}.

\bibitem[Nesterov(2015)]{Nesterov2015}
Y.~Nesterov.
\newblock Universal gradient methods for convex optimization problems.
\newblock \emph{Math. Program.}, 152\penalty0 (1-2):\penalty0 381--404, 2015.

\bibitem[Peng and Wei(2007)]{Peng2007}
J.~Peng and Y.~Wei.
\newblock Approximating {K}--means--type clustering via semidefinite
  programming.
\newblock \emph{SIAM J. Optim.}, 18\penalty0 (1):\penalty0 186--205, 2007.

\bibitem[Ryu and Boyd(2016)]{Ryu2016}
E.~K. Ryu and S.~Boyd.
\newblock Primer on monotone operator methods.
\newblock \emph{Appl. Comput. Math.}, 15\penalty0 (1):\penalty0 3--43, 2016.

\bibitem[Tran-Dinh et~al.(2018)Tran-Dinh, Fercoq, and Cevher]{TranDinh2017}
Q.~Tran-Dinh, O.~Fercoq, and V.~Cevher.
\newblock A smooth primal-dual optimization framework for nonsmooth composite
  convex minimization.
\newblock \emph{SIAM J. Optim.}, 28\penalty0 (1):\penalty0 96--134, 2018.

\bibitem[Wright(1997)]{Wright1997}
S.~J. Wright.
\newblock \emph{Primal--Dual Interior--Point Methods}.
\newblock SIAM, Philadelphia, USA, 1997.

\bibitem[Ye()]{GSet}
Y.~Ye.
\newblock Gset random graphs,
  https://www.cise.ufl.edu/research/sparse/matrices/gset/.

\bibitem[Yen et~al.(2016{\natexlab{a}})Yen, Huang, Zhong, Ravikumar, and
  Dhillon]{Yen2016svm}
I.~E.-H. Yen, K.~Huang, R.~Zhong, P.~Ravikumar, and I.~S. Dhillon.
\newblock Dual decomposed learning with factorwise oracle for structural svm
  with large output domain.
\newblock In \emph{Advances in Neural Information Processing Systems 29},
  2016{\natexlab{a}}.

\bibitem[Yen et~al.(2016{\natexlab{b}})Yen, Lin, Zhang, Ravikumar, and
  Dhillon]{yen2016convex}
I.~E.-H. Yen, X.~Lin, J.~Zhang, P.~Ravikumar, and I.~S. Dhillon.
\newblock A convex atomic--norm approach to multiple sequence alignment and
  motif discovery.
\newblock In \emph{Proceedings of the $33$rd International Conference on
  Machine Learning}, 2016{\natexlab{b}}.

\bibitem[Yurtsever et~al.(2015)Yurtsever, Tran-Dinh, and Cevher]{Yurtsever2015}
A.~Yurtsever, Q.~Tran-Dinh, and V.~Cevher.
\newblock A universal primal-dual convex optimization framework.
\newblock In \emph{Advances in Neural Information Processing Systems 28}, 2015.

\bibitem[Yurtsever et~al.(2017)Yurtsever, Udell, Tropp, and
  Cevher]{Yurtsever2016}
A.~Yurtsever, M.~Udell, J.~Tropp, and V.~Cevher.
\newblock Sketchy decisions: Convex low-rank matrix optimization with optimal
  storage.
\newblock In \emph{Proc. 20th Int. Conf. Artificial Intelligence and Statistics
  (AISTATS)}, 2017.

\bibitem[Yurtsever et~al.(2018)Yurtsever, Fercoq, Locatello, and
  Cevher]{Yurtsever2018}
A.~Yurtsever, O.~Fercoq, F.~Locatello, and V.~Cevher.
\newblock A conditional gradient framework for composite convex minimization
  with applications to semidefinite programming.
\newblock In \emph{Proc. $35$th Int. Conf. Machine Learning}, 2018.

\end{thebibliography}
\bibliographystyle{abbrvnat}

\newpage

\onecolumn
\appendix

\section{Proof of convergence}

For notational simplicity in the proof, we redefine augmented Lagrangian function with three variables, including the slack variable $r = Ax \in \mathcal{K}$ as
\begin{align*}\label{eqn:augmented-lagrangian-general}
\mathcal{L}_{\lambda}(x,r,y) :&= f(x) + \ip{y}{Ax - r} + \frac{\lambda}{2}\norm{Ax - r}^2,
\end{align*}
where $y \in \R^d$ is the Lagrange multiplier and $\lambda >0$ is the augmented Lagrangian parameter. 

Directional derivatives of augmented Lagrangian function can be written as
\begin{align*}
\nabla_x \mathcal{L}_{\lambda}(x,r,y) & = \nabla f(x)  + A^\top y + \lambda A^\top(Ax - r) \\
\nabla_r \mathcal{L}_{\lambda}(x,r,y) & = -y - \lambda(Ax-r)\\
\nabla_y \mathcal{L}_{\lambda}(x,r,y) & = Ax - r.
\end{align*}

Denote by $\bar{L}_k = (L_f + \lambda_k\norm{A}^2)$. Then, using the Taylor expansion, we get the following estimate:
\begin{align*}
\mathcal{L}_{\lambda_{k+1}}(x_{k+1},r_k,y_k) 
& \leq \mathcal{L}_{\lambda_{k+1}}(x_k,r_k,y_k) + \ip{\nabla_x \mathcal{L}_{\lambda_{k+1}}(x_k,r_k,y_k) }{x_{k+1} - x_k} + \frac{\bar{L}_{k+1}}{2} \norm{x_{k+1} - x_k }^2 \\
& = \mathcal{L}_{\lambda_{k+1}}(x_k,r_k,y_k) + \eta_k \ip{\nabla_x \mathcal{L}_{\lambda_{k+1}}(x_k,r_k,y_k) }{s_k - x_k} + \eta_k^2 \frac{\bar{L}_{k+1}}{2} \norm{s_k - x_k }^2 \\
& = \mathcal{L}_{\lambda_k}(x_k,r_k,y_k) + \eta_k \ip{\nabla_x \mathcal{L}_{\lambda_k}(x_k,r_k,y_k) }{s_k - x_k} + \eta_k^2 \frac{\bar{L}_{k+1}}{2} D_{\mathcal{X}}^2\\
& \qquad + \frac{\lambda_{k+1} -\lambda_k}{2}\norm{Ax_k-r_k}^2 + \eta_k (\lambda_{k+1} - \lambda_k)\ip{Ax_k - r_k}{As_k-Ax_k}.
\end{align*}

We can bound the inner product term on the right hand as follows:
\begin{align*}
\ip{\nabla_x \mathcal{L}_{\lambda_k}(x_k,r_k,y_k)}{s_k - x_k} 
& = \ip{\nabla f(x_k) + \lambda_k A^\top (Ax_k - r_k) + A^\top y_k}{s_k - x_k} \\
& \leq \ip{\nabla f(x_k) + \lambda_k A^\top (Ax_k - r_k) +  A^\top y_k}{x^\star - x_k} \\
& = \ip{\nabla f(x_k)}{x^\star - x_k} + \ip{\lambda_k (Ax_k - r_k) + y_k}{A x^\star - A x_k} \\
& = \ip{\nabla f(x_k)}{x^\star - x_k} + \ip{\lambda_k (Ax_k - r_k) + y_k}{r_k - r_k + A x^\star - A x_k} \\
& = \ip{\nabla f(x_k)}{x^\star - x_k} -  \lambda_k\norm{ Ax_k - r_k}^2 - \ip{y_k}{Ax_k - r_k} \\
&\qquad\qquad+ \lambda_k \ip{Ax_k-r_k}{Ax^\star - r_k} + \ip{y_k}{Ax^\star - r_k} \\
& \leq f^\star - f(x_k) -  \lambda_k\norm{ Ax_k - r_k}^2 - \ip{y_k}{Ax_k - r_k} \\
&\qquad\qquad+ \lambda_k \ip{Ax_k-r_k}{Ax^\star - r_k} + \ip{y_k}{Ax^\star - r_k} \\
& = f^\star - \mathcal{L}_{\lambda_k} (x_k,r_k,y_k) - \frac{\lambda_k}{2}\norm{ Ax_k - r_k}^2 + \lambda_k \ip{Ax_k-r}{Ax^\star - r_k} + \ip{y_k}{Ax^\star - r_k} \\
& = \mathcal{L}^\star - \mathcal{L}_{\lambda_k} (x_k,r_k,y_k) + \frac{\lambda_k}{2} \norm{Ax^\star - r_k}^2 - \frac{\lambda_k}{2}\norm{Ax_k - Ax^\star}^2 + \ip{y_k}{Ax^\star - r_k} .
\end{align*}
where the first inequality holds since $s_k$ is the solution of $\mathrm{lmo}$, the second inequality simply follows the convexity of $f$, and the last equality holds due to strong duality. 

Also note by definition, $\bar{r}_{k+1} = \arg\min_{r \in \mathcal{K}} \mathcal{L}_{\lambda_{k+1}}(x_{k+1},r,y_k)$, hence
\begin{align*}
\mathcal{L}_{\lambda_{k+1}}(x_{k+1},\bar{r}_{k+1},y_k)
\leq 
\mathcal{L}_{\lambda_{k+1}}(x_{k+1},r_k,y_k).
\end{align*}

Combining these bounds, we arrive at
\begin{align*}
\mathcal{L}_{\lambda_{k+1}}(x_{k+1},\bar{r}_{k+1},y_k) 
& \leq (1-\eta_k) \mathcal{L}_{\lambda_k}(x_k,r_k,y_k)  + (1-\eta_k) \left(\frac{\lambda_{k+1}}{2} - \frac{\lambda_k}{2}\right) \norm{Ax_k-r_k}^2 + \eta_k \mathcal{L}^\star + \eta_k^2 \frac{\bar{L}_{k+1}}{2}  D_{\mathcal{X}}^2  \\
&\qquad \qquad+\eta_k\frac{\lambda_k}{2} \norm{Ax^\star - r}^2 + \eta_k\ip{y_k}{Ax^\star - r_k} - \eta_k\frac{\lambda_k}{2} \norm{Ax_k - Ax^\star}^2 \\
& \leq (1-\eta_k) \mathcal{L}_{\lambda_k}(x_k,r_k,y_k)  + (1-\eta_k) \left(\frac{\lambda_{k+1}}{2} - \frac{\lambda_k}{2}\right) \norm{Ax_k-r_k}^2 + \eta_k \mathcal{L}^\star + \eta_k^2 \frac{\bar{L}_{k+1}}{2}  D_{\mathcal{X}}^2  \\
&\qquad \qquad- \eta_k\frac{\lambda_k}{2} \norm{Ax_k - r_k}^2 + \eta_k\ip{y_k + \lambda_k(Ax_k - r_k)}{Ax^\star - r_k}\\
& \leq (1-\eta_k) \mathcal{L}_{\lambda_k}(x_k,r_k,y_k)  + (1-\eta_k) \left(\frac{\lambda_{k+1}}{2} - \frac{\lambda_k}{2}\right) \norm{Ax_k-r_k}^2 + \eta_k \mathcal{L}^\star + \eta_k^2 \frac{\bar{L}_{k+1}}{2}  D_{\mathcal{X}}^2  \\
&\qquad \qquad- \eta_k\frac{\lambda_k}{2} \norm{Ax_k - r_k}^2 
\end{align*}
where the last inequality follows from the optimality condition. By definition, $r_k = \arg\min_{r \in \mathcal{K}} \mathcal{L}_{\lambda_k}(x_k,r,y_k)$, hence the following estimate holds $\forall r \in \mathcal{K}$
\begin{align*}
\ip{y_k + \lambda_k(Ax_k - r_k)}{r - r_k} = \ip{\nabla_r \mathcal{L}_{\lambda_k}(x_k,r_k,y_k)}{r - r_k} \leq 0
\end{align*}
and in particular for $r = Ax^\star \in \mathcal{K}$. 

In order to obtain a recurrence, we need to shift the dual variable on the left hand side of our bound. For this, we use the following relations:
\begin{align*}
\mathcal{L}_{\lambda_{k+1}}(x_{k+1},r_{k+1},y_{k+1}) 
& \leq \mathcal{L}_{\lambda_{k+1}}(x_{k+1},\bar{r}_{k+1},y_{k+1})  \\
& = \mathcal{L}_{\lambda_{k+1}}(x_{k+1},\bar{r}_{k+1},y_{k}) + \ip{y_{k+1} - y_k}{Ax_{k+1} - \bar{r}_{k+1}}  \\
& = \mathcal{L}_{\lambda_{k+1}}(x_{k+1},\bar{r}_{k+1},y_{k}) + \sigma_{k+1} \norm{Ax_{k+1} - \bar{r}_{k+1}}^2.
\end{align*}
Combining all these bounds and subtracting $\mathcal{L}^\star$ from both sides, we end up with
\begin{equation}\label{eqn:proof-split-point}
\begin{aligned} 
\mathcal{L}_{\lambda_{k+1}}(x_{k+1},r_{k+1},y_{k+1}) - \mathcal{L}^\star
& \leq (1-\eta_k) \left( \mathcal{L}_{\lambda_k}(x_k,r_k,y_k) - \mathcal{L}^\star\right)  + \eta_k^2 \frac{\bar{L}_{k+1}}{2}  D_{\mathcal{X}}^2   \\
&\qquad \qquad  + \frac{1}{2} \left((1-\eta_k)\lambda_{k+1} - \lambda_k \right) \norm{Ax_k-r_k}^2 + \sigma_{k+1} \norm{Ax_{k+1} - \bar{r}_{k+1}}^2
\end{aligned}
\end{equation}

From this point, we consider two cases: constant step size with growth condition, and decreasing step size. 

\subsection{Constant bound on step-size}
 
We choose $\sigma_{k+1} \leq \lambda_0$ by line-search, to ensure the following conditions:
\begin{align*}
\sigma_{k+1} \norm{Ax_{k+1} - \bar{r}_{k+1}}^2 \leq \eta_k^2 \frac{\bar{L}_{k+1}}{2}  D_{\mathcal{X}}^2 
\qquad \& \qquad
\norm{y_{k+1}} \leq D_{\mathcal{Y}_{k+1}}
\end{align*}
where $D_{\mathcal{Y}_{k+1}}$ is a sequence of positive and non-decreasing numbers, to input. Note that $\sigma_{k+1}$ is well defined, in the sense there exists $\sigma_{k+1} \geq 0$ which satisfy both conditions, simply because $\sigma_{k+1} = 0$ trivially satisfies them. 

In addition, since we choose $\lambda_k = \lambda_0 \sqrt{k+1} $ and $\eta_k = 2/{k+1}$, we have
\begin{align*}
(1-\eta_k)\lambda_{k+1} - \lambda_k = \frac{k-1}{k+1} \sqrt{k+2} - \sqrt{k+1} \leq \frac{k}{\sqrt{k+2}} - \sqrt{k+1} \leq 0.
\end{align*}

As a consequence, we can simplify \eqref{eqn:proof-split-point} as
\begin{align*}
\mathcal{L}_{\lambda_{k+1}}(x_{k+1},r_{k+1},y_{k+1}) - \mathcal{L}^\star
& \leq (1-\eta_k) \left( \mathcal{L}_{\lambda_k}(x_k,r_k,y_k) - \mathcal{L}^\star\right)  + \eta_k^2 \bar{L}_{k+1}  D_{\mathcal{X}}^2  
\end{align*}

Applying this recursion we get
\begin{align*}
\mathcal{L}_{\lambda_{k+1}}(x_{k+1},r_{k+1},y_{k+1})  - \mathcal{L}^\star
& \leq \prod_{j=1}^{k} (1-\eta_j) \left(\mathcal{L}_{\lambda_{1}}(x_1,r_1,y_1) - \mathcal{L}^\star\right) +  D_{\mathcal{X}}^2 \sum_{\ell=1}^{k} \eta_{\ell}^2 \bar{L}_{\ell+1}  \prod_{j=\ell}^k (1-\eta_j) \\
& = D_{\mathcal{X}}^2 \sum_{\ell=1}^{k} \eta_{\ell}^2 \bar{L}_{\ell+1}  \prod_{j=\ell}^k (1-\eta_j) 
\end{align*}
where the last equality follows since $\eta_1 = 1$. By using the following inequality,
\begin{align*}
\sum_{\ell=1}^{k} \eta_{\ell}^2 \prod_{j=\ell}^k (1-\eta_j) = \sum_{\ell=1}^k \frac{4}{(\ell + 1)^2} \prod_{j=\ell}^k \frac{j-1}{j+1} = \sum_{\ell=1}^k \frac{4}{(\ell+1)^2} \frac{(\ell-1)\ell}{k(k+1)} \leq \frac{4}{k+1},
\end{align*} 
we get the following bound on the augmented Lagrangian:
\begin{align*}
\mathcal{L}_{\lambda_{k+1}}(x_{k+1},r_{k+1},y_{k+1})  - \mathcal{L}^\star 
 \leq \frac{4}{k+1} D_{\mathcal{X}}^2 \bar{L}_{k+1} 
 =  4 D_{\mathcal{X}}^2 \left( \frac{L_f}{k+1}   + \frac{\lambda_0 \norm{A}^2}{\sqrt{k+1}} \right).
\end{align*}
In the next step, we translate the bound on augmented Lagrangian to convergence guarantees on objective residual and feasibility gap.

\paragraph{Convergence of objective.}
We start by using the definition of augmented Lagrangian and the strong duality:
\begin{align*}
f(x_{k+1}) - f^\star & = \mathcal{L}_{\lambda_{k+1}}(x_{k+1},r_{k+1},y_{k+1}) - \mathcal{L}^\star + \frac{1}{2\lambda_{k+1}}\norm{y_{k+1}}^2 -  \frac{\lambda_{k+1}}{2}\dist^2\Big(Ax_{k+1} + \frac{1}{\lambda}y_{k+1},\mathcal{K}\Big) \\
& \leq \mathcal{L}_{\lambda_{k+1}}(x_{k+1},r_{k+1},y_{k+1}) - \mathcal{L}^\star + \frac{D^2_{\mathcal{Y}_{k+1}}}{2\lambda_{k+1}}\\
& \leq 4 D_{\mathcal{X}}^2 \left( \frac{L_f}{k+1}   + \frac{\lambda_0 \norm{A}^2}{\sqrt{k+1}} \right) + \frac{D^2_{\mathcal{Y}_{k+1}}}{2\lambda_0\sqrt{k+1}}.
\end{align*}

For the lower bound, we use the classical Lagrange saddle point properties, that $\forall (x,r) \in \mathcal{X} \times \mathcal{K}$ we have
\begin{align}\label{eqn:lagrange-saddle-property}
f^\star \leq \mathcal{L}(x,r,y^\star) = f(x) + \ip{y^\star}{Ax - r} \leq f(x) + \norm{y^\star}\norm{Ax - r}.
\end{align}
By choosing $x = x_{k+1}$ and $r = \proj{\mathcal{K}}(Ax_{k+1})$ and rearranging, we arrive at
\begin{align*}
f(x_{k+1}) - f^\star \geq - \norm{y^\star} \dist \left(Ax_{k+1} ,\mathcal{K} \right).
\end{align*}

\paragraph{Convergence of feasibility.} We start by combining \eqref{eqn:proof-split-point} and \eqref{eqn:lagrange-saddle-property} by choosing $x = x_{k+1}$ and $r = r_{k+1}$:
\begin{gather*}
 \ip{y_{k+1} - y^\star}{Ax_{k+1}-r_{k+1}} + \frac{\lambda_{k+1}}{2} \norm{Ax_{k+1}-r_{k+1}}^2 \leq \frac{4}{k+1} D_{\mathcal{X}}^2 \bar{L}_{k+1} \\
\implies -\norm{y_{k+1} - y^\star}\norm{Ax_{k+1}-r_{k+1}} + \frac{\lambda_{k+1}}{2} \norm{Ax_{k+1}-r_{k+1}}^2 \leq \frac{4}{k+1} D_{\mathcal{X}}^2 \bar{L}_{k+1} 
\end{gather*}
This is a second order inequality with respect to $\norm{Ax_{k+1}-r_{k+1}}$, and by solving this inequality we get
\begin{align*}
\norm{Ax_{k+1}-r_{k+1}}
& \leq \frac{1}{\lambda_{k+1}} \left( \norm{y_{k+1} - y^\star} + \sqrt{\norm{y_{k+1} - y^\star}^2 + \frac{8D^2_{\mathcal{X}}\bar{L}_{k+1}}{k+1}\lambda_{k+1}}  \right) \\
& \leq \frac{1}{\lambda_{k+1}} \left( 2\norm{y_{k+1} - y^\star} + \sqrt{\frac{8 D^2_{\mathcal{X}}}{k+1} (L_f + \lambda_{k+1})\lambda_{k+1}}  \right) \\
& \leq \frac{2}{\lambda_0\sqrt{k+1}} \left( \norm{y_{k+1} - y^\star} + D_{\mathcal{X}}\sqrt{2 \left(\frac{L_f \lambda_0}{\sqrt{k+1}} + \norm{A}^2\lambda_0^2\right) } \right).
\end{align*}

Finally, we use the properties of projection to get the bound on the feasibility gap:
\begin{gather*}
\dist \left(Ax_{k+1},\mathcal{K} \right) = \norm{Ax_{k+1} - \dot{r}_{k+1}} = \norm{Ax_{k+1} - r_{k+1} + r_{k+1} - \dot{r}_{k+1}} \leq \norm{Ax_{k+1} - r_{k+1}} + \norm{r_{k+1} - \dot{r}_{k+1}} \\
 \leq \norm{Ax_{k+1} - r_{k+1}} + \norm{Ax_{k+1} - Ax_{k+1} + \frac{1}{\lambda_{k+1}}y_{k+1}} \leq \norm{Ax_{k+1} - r_{k+1}} + \frac{D_{\mathcal{Y}_{k+1}}}{\lambda_{k+1}}.
\end{gather*}

\subsection{Decreasing bound on step-size}

Choose parameters
\begin{align*}
\lambda_k = \lambda_0 \sqrt{k+1} \qquad \sigma_k = \frac{\lambda_0}{2\sqrt{k+1}} \qquad \eta_k = \frac{2}{k+1}.
\end{align*}

Now we execute the last term using the non-expansiveness of projection operator
\begin{align*}
\norm{Ax_{k+1} - \bar{r}_{k+1}}
& = \norm{r_{k+1} - r_{k+1}  + Ax_{k+1} - \bar{r}_{k+1}} \\
& \leq \norm{r_{k+1} - \bar{r}_{k+1}}  + \norm{Ax_{k+1} - r_{k+1}}\\
& \leq \frac{1}{\lambda_{k+1}}\norm{y_{k+1} - y_k}  + \norm{Ax_{k+1} - r_{k+1}}\\
& = \frac{\sigma_{k+1}}{\lambda_{k+1}}\norm{Ax_{k+1} - \bar{r}_{k+1}}  + \norm{Ax_{k+1} - r_{k+1}},
\end{align*}
hence $\norm{Ax_{k+1} - \bar{r}_{k+1}} \leq (1-\sigma_{k+1}/\lambda_{k+1})^{-1}\norm{Ax_{k+1} - r_{k+1}}$.

Overall, we obtain the following recursion relation:
\begin{align*}
\mathcal{L}_{\lambda_{k+1}}(x_{k+1},r_{k+1},y_{k+1}) - \mathcal{L}^\star
& \leq (1-\eta_k) \left( \mathcal{L}_{\lambda_k}(x_k,r_k,y_k) - \mathcal{L}^\star\right)  + \eta_k^2 \frac{\bar{L}_{k+1}}{2}  D_{\mathcal{X}}^2   \\
&\qquad \qquad  + \frac{1}{2} \left((1-\eta_k)\lambda_{k+1} - \lambda_k \right) \norm{Ax_k-r_k}^2 + \frac{\lambda_{k+1}^2\sigma_{k+1}}{(\lambda_{k+1}-\sigma_{k+1})^2} \norm{Ax_{k+1} - r_{k+1}}^2
\end{align*}
Now we can apply recursion, and we get
\begin{align*}
\mathcal{L}_{\lambda_{k+1}}(x_{k+1},r_{k+1},y_{k+1})  - \mathcal{L}^\star
& \leq \prod_{j=1}^{k} (1-\eta_j) \left(\mathcal{L}_{\lambda_{1}}(x_1,r_1,y_1) - \mathcal{L}^\star\right) + \frac{ D_{\mathcal{X}}^2}{2} \sum_{\ell=1}^{k} \eta_{\ell}^2 \bar{L}_{\ell+1}  \prod_{j=\ell}^k (1-\eta_j) \\
& \qquad\qquad + \frac{1}{2} \sum_{\ell=1}^{k} \left((1-\eta_{\ell})\lambda_{\ell+1} - \lambda_{\ell} \right) \norm{Ax_{\ell}-r_{\ell}}^2 \prod_{j=\ell}^k (1-\eta_j) \\
& \qquad\qquad + \sum_{\ell=1}^{k} \frac{\lambda^2_{\ell+1}\sigma_{\ell+1}}{(\lambda_{\ell+1}-\sigma_{\ell+1})^2} \norm{Ax_{\ell+1} - r_{\ell+1}}^2 \prod_{j=\ell}^k (1-\eta_j) .
\end{align*}
Note that the terms which involve $(1-\eta_1)$ on the right hand side are zero since $\eta_1 = 1$. 

Now, we focus on the last summation term:
\begin{align*}
\sum_{\ell=1}^{k} \frac{\lambda^2_{\ell+1}\sigma_{\ell+1}}{(\lambda_{\ell+1}-\sigma_{\ell+1})^2}&  \norm{Ax_{\ell+1} - r_{\ell+1}}^2 \prod_{j=\ell}^k (1-\eta_j)  \\
& =  \sum_{\ell=1}^{k} (1-\eta_{\ell}) \frac{\lambda^2_{\ell+1}\sigma_{\ell+1}}{(\lambda_{\ell+1}-\sigma_{\ell+1})^2} \norm{Ax_{\ell+1} - r_{\ell+1}}^2 \prod_{j=\ell+1}^k (1-\eta_j) \\
&  = \sum_{\ell=2}^{k+1} (1-\eta_{\ell-1}) \frac{\lambda^2_{\ell}\sigma_{\ell}}{(\lambda_{\ell}-\sigma_{\ell})^2} \norm{Ax_{\ell} - r_{\ell}}^2 \prod_{j=\ell}^{k+1} (1-\eta_j) \\
&  \leq \sum_{\ell=2}^{k} (1-\eta_{\ell-1}) \frac{\lambda^2_{\ell}\sigma_{\ell}}{(\lambda_{\ell}-\sigma_{\ell})^2} \norm{Ax_{\ell} - r_{\ell}}^2 \prod_{j=\ell}^k (1-\eta_j) + (1-\eta_k) \frac{\lambda^2_{k+1}\sigma_{k+1}}{(\lambda_{k+1} - \sigma_{k+1})^2} \norm{Ax_{k+1} - r_{k+1}}^2.
\end{align*}

We choose parameters $\lambda_k, \eta_k$ and $\sigma_{k+1}$ such that for all $k \geq 2$, we have 
\begin{align*}
\bigg( \frac{1}{2} \left((1-\eta_{\ell})\lambda_{\ell+1} - \lambda_{\ell} \right) + (1-\eta_{\ell-1})\frac{\lambda^2_{\ell}\sigma_{\ell}}{(\lambda_{\ell}-\sigma_{\ell})^2} \bigg) \leq 0,
\end{align*}
hence by combining these bounds, we get 
\begin{align*}
\mathcal{L}_{\lambda_{k+1}}(x_{k+1},r_{k+1},y_{k+1})  - \mathcal{L}^\star \leq \frac{ D_{\mathcal{X}}^2}{2} \sum_{\ell=1}^{k} \eta_{\ell}^2 \bar{L}_{\ell+1}  \prod_{j=\ell}^k (1-\eta_j)  +(1-\eta_k) \frac{\lambda^2_{k+1}\sigma_{k+1}}{(\lambda_{k+1} - \sigma_{k+1})^2} \norm{Ax_{k+1} - r_{k+1}}^2  .
\end{align*}

Using the following formula
\begin{align*}
\sum_{\ell=1}^{k} \eta_{\ell}^2 \prod_{j=\ell}^k (1-\eta_j) = \sum_{\ell=1}^k \frac{4}{(\ell + 1)^2} \prod_{j=\ell}^k \frac{j-1}{j+1} = \sum_{\ell=1}^k \frac{4}{(\ell+1)^2} \frac{(\ell-1)\ell}{k(k+1)} \leq \frac{4}{k+1}
\end{align*} 
we get the following bound on the augmented Lagrangian:
\begin{align*}
\mathcal{L}_{\lambda_{k+1}}(x_{k+1},r_{k+1},y_{k+1})  - \mathcal{L}^\star 
& \leq \frac{2}{k+1} D_{\mathcal{X}}^2 \bar{L}_{k+1} +(1-\eta_k) \frac{\lambda^2_{k+1}\sigma_{k+1}}{(\lambda_{k+1} - \sigma_{k+1})^2} \norm{Ax_{k+1} - r_{k+1}}^2\\
& =  2 D_{\mathcal{X}}^2 \left( \frac{L_f}{k+1}   + \frac{\lambda_0 \norm{A}^2}{\sqrt{k+1}} \right) + \frac{\lambda_0}{2\sqrt{k+1}} \norm{Ax_{k+1} - r_{k+1}}^2.
\end{align*}

\paragraph{Convergence of objective.} Lower bound of the objective residual follows similarly to the constant step-size case. 
For upper bound, we start by
\begin{align*}
f(x_{k+1}) - f^\star & = \mathcal{L}_{\lambda_{k+1}}(x_{k+1},r_{k+1},y_{k+1}) - \mathcal{L}^\star - \ip{y_{k+1}}{Ax_{k+1} - r_{k+1}} - \frac{\lambda_{k+1}}{2} \norm{Ax_{k+1} - r_{k+1}}^2\\
& \leq 2 D_{\mathcal{X}}^2 \left( \frac{L_f}{k+1}   + \frac{\lambda_0 \norm{A}^2}{\sqrt{k+1}} \right) + \frac{\lambda_0}{2} \left( \frac{1}{\sqrt{k+1}} - \sqrt{k+1} \right) \norm{Ax_{k+1} - r_{k+1}}^2 - \ip{y_{k+1}}{Ax_{k+1} - r_{k+1}} \\
& \leq 2 D_{\mathcal{X}}^2 \left( \frac{L_f}{k+1}   + \frac{\lambda_0 \norm{A}^2}{\sqrt{k+1}} \right) - \frac{\lambda_k}{2} \norm{Ax_{k+1} - r_{k+1}}^2 + \norm{y_{k+1}}\norm{Ax_{k+1} - r_{k+1}} \\
& \leq 2 D_{\mathcal{X}}^2 \left( \frac{L_f}{k+1}   + \frac{\lambda_0 \norm{A}^2}{\sqrt{k+1}} \right) + \frac{\norm{y_{k+1}}^2}{\lambda_0 \sqrt{k+1}}
\end{align*}

\paragraph{Convergence of feasibility.} We start by combining \eqref{eqn:proof-split-point} and \eqref{eqn:lagrange-saddle-property} by choosing $x = x_{k+1}$ and $r = r_{k+1}$:
\begin{gather*}
 \ip{y_{k+1} - y^\star}{Ax_{k+1}-r_{k+1}} + \frac{\lambda_{k+1}}{2} \norm{Ax_{k+1}-r_{k+1}}^2 \leq \frac{2}{k+1} D_{\mathcal{X}}^2 \bar{L}_{k+1} +\frac{\lambda_0}{2\sqrt{k+1}} \norm{Ax_{k+1} - r_{k+1}}^2 \\
\implies -\norm{y_{k+1} - y^\star}\norm{Ax_{k+1}-r_{k+1}} + \frac{\lambda_{k}}{2} \norm{Ax_{k+1}-r_{k+1}}^2 \leq \frac{2}{k+1} D_{\mathcal{X}}^2 \bar{L}_{k+1} 
\end{gather*}
This is a second order inequality with respect to $\norm{Ax_{k+1}-r_{k+1}}$, and by solving this inequality we get
\begin{align*}
\norm{Ax_{k+1}-r_{k+1}}
& \leq \frac{1}{\lambda_{k}} \left( \norm{y_{k+1} - y^\star} + \sqrt{\norm{y_{k+1} - y^\star}^2 + \frac{4D^2_{\mathcal{X}}\bar{L}_{k+1}}{k+1}\lambda_{k}}  \right) \\
& \leq \frac{2}{\lambda_{k}} \left( \norm{y_{k+1} - y^\star} + \sqrt{\frac{D^2_{\mathcal{X}}}{k+1} (L_f + \lambda_{k+1})\lambda_{k}}  \right) \\
& \leq \frac{2}{\lambda_0\sqrt{k+1}} \left( \norm{y_{k+1} - y^\star} + D_{\mathcal{X}}\sqrt{\frac{L_f}{k+1} + \lambda_0^2 } \right).
\end{align*}

To complete the proof, we use the following arguments: Denote by $\dot{r}_{k+1} = \arg\min_{r\in\mathcal{K}} \norm{Ax_{k+1} - r}$, we have
\begin{gather*}
\dist \left(Ax_{k+1},\mathcal{K} \right) = \norm{Ax_{k+1} - \dot{r}_{k+1}} = \norm{Ax_{k+1} - r_{k+1} + r_{k+1} - \dot{r}_{k+1}} \leq \norm{Ax_{k+1} - r_{k+1}} + \norm{r_{k+1} - \dot{r}_{k+1}} \\
 \leq \norm{Ax_{k+1} - r_{k+1}} + \norm{Ax_{k+1} - Ax_{k+1} + \frac{1}{\lambda_{k+1}}y_{k+1}} \leq \norm{Ax_{k+1} - r_{k+1}} + \frac{1}{\lambda_{k+1}}\norm{y_{k+1}}.
\end{gather*}

\newpage

\section{Additional Numerical Experiments \& Observations}

This appendix presents implementation details and additional results from the numerical experiments in \Cref{sec:experiments}. 

In the last two pages of this document, we present all trials of each method for max-cut SDP with GD97\_b dataset. Note that FW-AL and IAL has more $2$ parameters to tune. We denote by $\rho$ in the legends of FW-AL plots the ratio between $\sigma_0$ and $2/\lambda$. Note that the method requires $\sigma_0 \leq \min\{\frac{2}{\lambda}, \frac{\alpha^2}{2\delta} \}$ to be tuned, where $\alpha$ is unknown and can be arbitrarily small. 

We also provide a brief conclusion about all methods and our observations:

{\bf CGAL \eqref{eqn:dual-step-constant}} and {\bf CGAL \eqref{eqn:dual-step-decreasing}:} CGAL with \eqref{eqn:dual-step-decreasing} step variant outperforms the base method, HCGM in this experiment as well as other experiments we considered. Nevertheless, we did not encounter any instance where CGAL \eqref{eqn:dual-step-decreasing} outperforms CGAL \eqref{eqn:dual-step-constant}, hence we focus on the \eqref{eqn:dual-step-constant} step variant. Note however CGAL \eqref{eqn:dual-step-decreasing} is still interesting from a theoretical perspective. Similar to the $\mathcal{O}(1/\sqrt{k})$ convergence guarantee of the feasibility gap, we can also show that the norm of updates $\|Ax_{k+1} - \bar{r}_{k+1}\|$ is decreasing with $\mathcal{O}(1/\sqrt{k})$ rate. Coupled with decreasing step size of the same rate and by triangle inequality, we can bound the norm of $\|y_{k+1}\|$ as the sum of terms that we expect to decrease by $\mathcal{O}(1/k)$ rate, resulting in a logarithmic bound naturally, without further conditions. We also did not encounter any problem in CGAL (both cases) even when we completely remove the conditions on $\|y_{k+1}\|$ in various tests. Unfortunately, we do not have guarantees for this case for now. 

{\bf HCGM:} HCGM is the base method for CGAL, and can be recovered from CGAL simply by choosing $y1=0$ and $\sigma_k = 0$.
HCGM guarantees $\mathcal{O}(1/\sqrt{k})$ convergence rate in the objective residual and the feasibility gap, which is optimal according to \cite{Yurtsever2015}, in ithe sense it matches the best rate for smoothness of the Lafrange dual problem. HCGM is a very simple method, easy to analyze, interpret and tune, but as we observed in various numerical experiments, this method typically performs with the worst case bounds in practice. CGAL specifically focuses on the practical performance and implementation of HCGM, extending it from quadratic penalty to augmented Lagrangian setup. As a result, CGAL retains essentially the same guarantees as HCGM, but performs much better in practice, achieving $\mathcal{O}(1/k)$ empirical rate in most instances. 

{\bf FW-AL:} FW-AL iterations are similar to CGAL, but the penalty parameter is fixed in contrast. The method, hence directly relies on the Lagrange multiplier for the convergence. This requires strong assumptions such as error bounds, and the theoretical analysis of this method is much more complicated than CGAL. The bounds are non-adaptive and depends on the unknown error bound parameter $\alpha$, which is proved to be positive assuming that Slater's condition holds. Nevertheless, this unknown constant directly appears in the bounds and the parameters. We argue that this constant can be arbitrarily small, and this method might be not implementable in practice. Even for a small scale max-cut SDP problem, after extensive search of proper parameters, we failed to find good parameter choices for this method, supporting our arguments. 

{\bf IAL:} This method theoretically has $\mathcal{O}(1/\epsilon^4)$ complexity of lmo calls. Nevertheless, the method performs better in practice, but requires a lot of effort for tuning. The method has a double-loop structure, and only the outer iterates provide reasonable approximations (which results in the stair like plots). We also tried evaluating the performance of the inner iterations, but the method simply jumps at the beginning of each subproblem due to CGM initialization. We also tried line-search to avoid this, but the method performs worse with line-search overall.  

{\bf UPD} and {\bf AUPD:} Remark that the problem instances we consider have bounded subgradients in the simple Lagrange dual formulation due to boundedness of domain $\mathcal{X}$. This corresponds to $0$-th order H\"older smoothness in the dual, hence UPD and AUPD both have $\mathcal{O}(1/\sqrt{k})$ rate of convergence, which is optimal according to \cite{Yurtsever2015}. Important to underline once again, that these methods are proved to converge only up to some accuracy level $\epsilon$. Indeed, we can easily observe this saturation in the objective residual in various of our numerical experiments. 

In our numerical experiment with small max-cut dataset, we observed similar performance of UPD and AUPD in terms of convergence rate, which is expected since the dual is only $0$-th order H\"older smooth. Interestingly, saturation of AUPD is not observed in contrast with the guarantees. One simple explanation for this observation is as follows: Both UPD and AUPD uses an inexact line-search procedure, but UPD lets $\delta = \epsilon$ error at each iteration, while AUPD requires increasing accuracy in line-search and only lets $\delta \sim \epsilon/k$ error in $k$th iteration. This decrease is required from the theoretical point of view to prevent error accumulation due to acceleration, but might not occur in practice, at least until later iterations. When error accumulation does not occur, decreasing inexactness also prevents saturation of UPD. Note however this comes at an increased computational cost. Dual objective depends on the output of sharp operator, and errors in sharp operator directly translates as objective evaluation. Considering the decreasing amount of inexactness, this method requires very accurate evaluations of the sharp operator. This does not cause much problem in very small scale problems where we can compute lmo in the exact sense, but even in medium scale problems with $1000 \times 1000$ dimensions we observed AUPD getting stucked in the line-search condition. Note that when the error in dual objective evaluations in two consecutive iterations is larger than the inexactness parameter $\delta$, line-search condition may become ill-defined, in the sense line-search turns into an infinite loop. 
Since UPD performs similarly as AUPD in small scale experiments and due to its robustness compared to AUPD, we focus on UPD for other experiments. 

\begin{figure}[ht]
\includegraphics[width=0.47\textwidth]{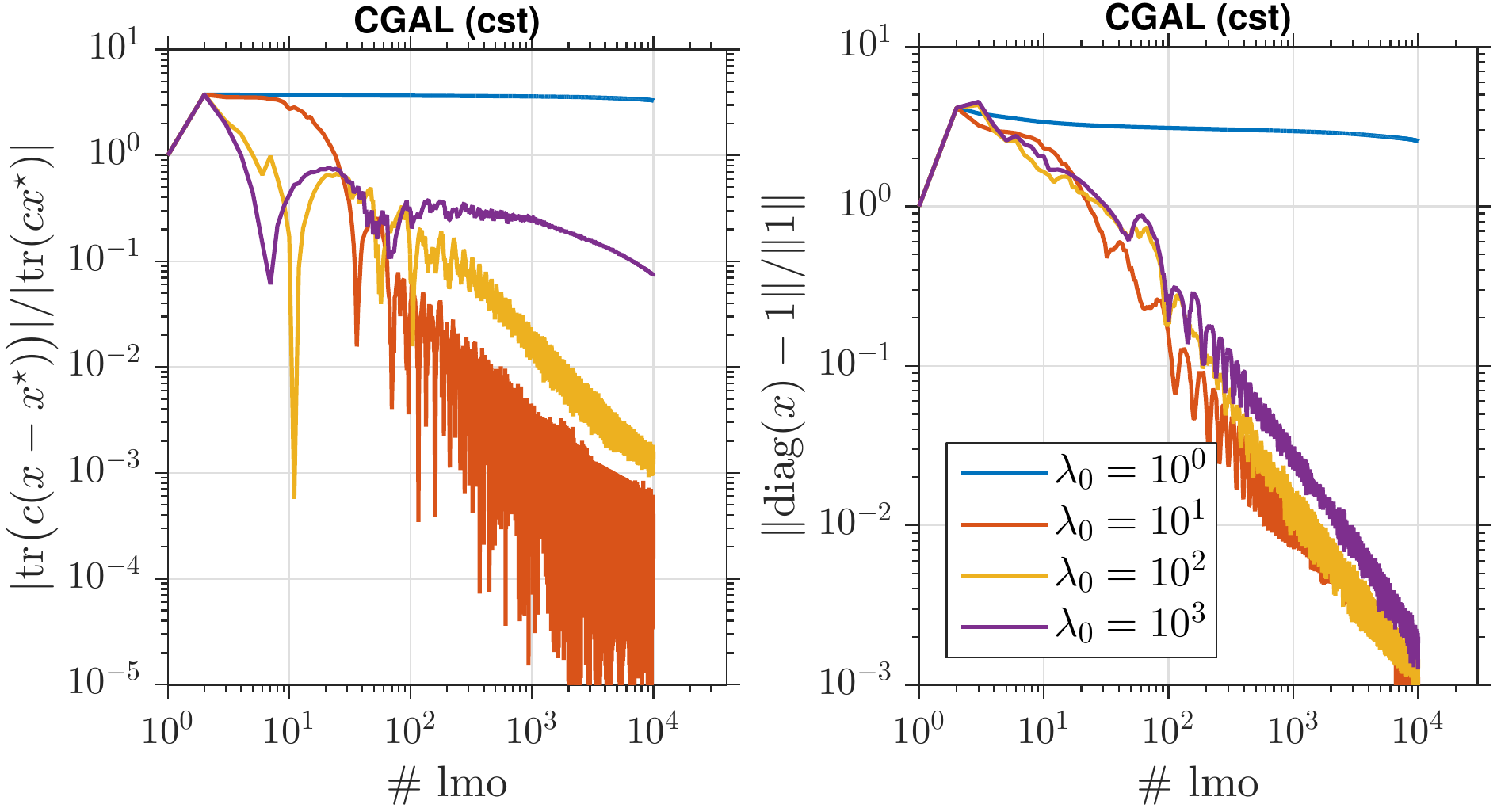}~~~~
\includegraphics[width=0.47\textwidth]{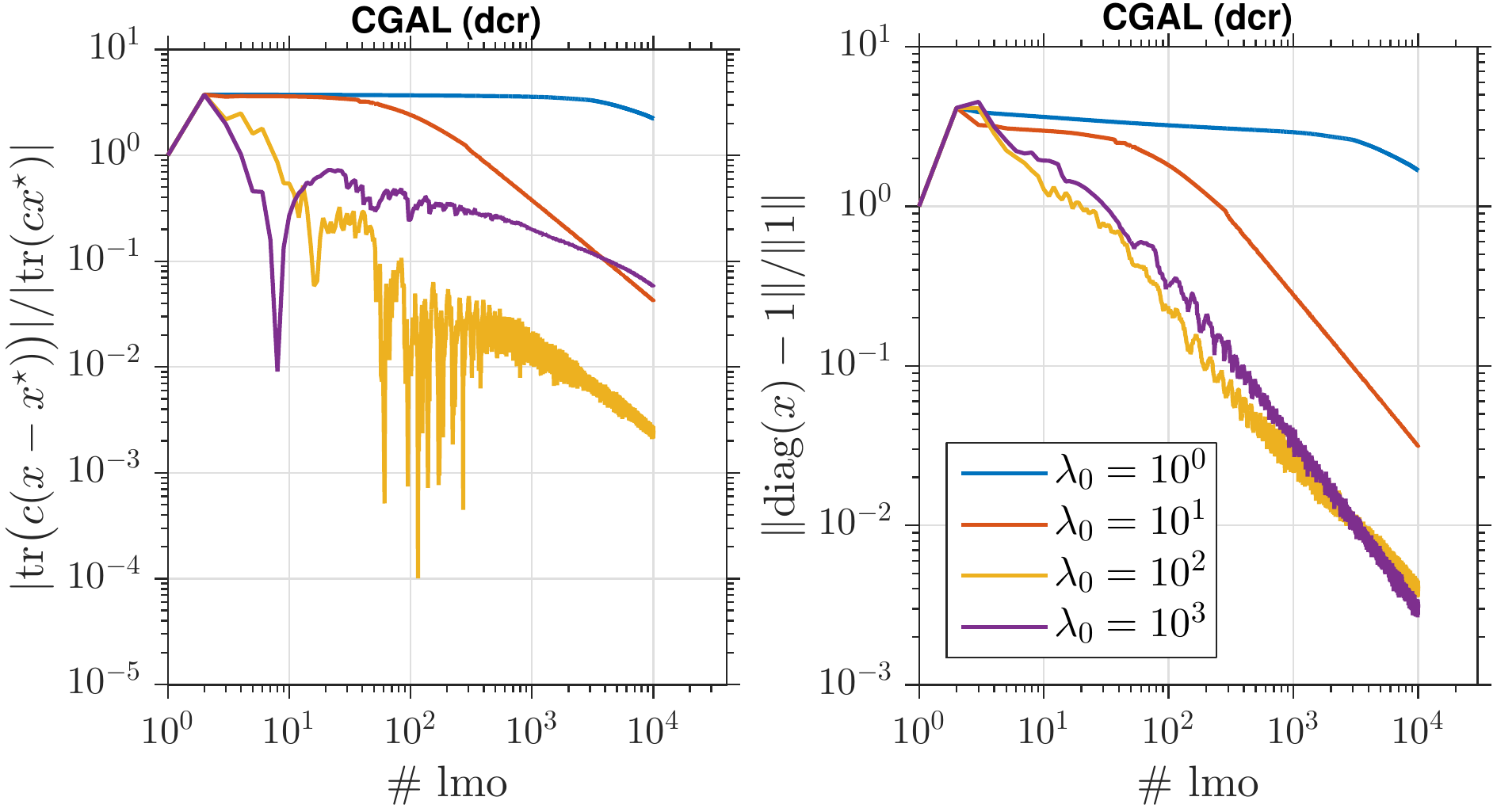}\\[1em]
\includegraphics[width=0.47\textwidth]{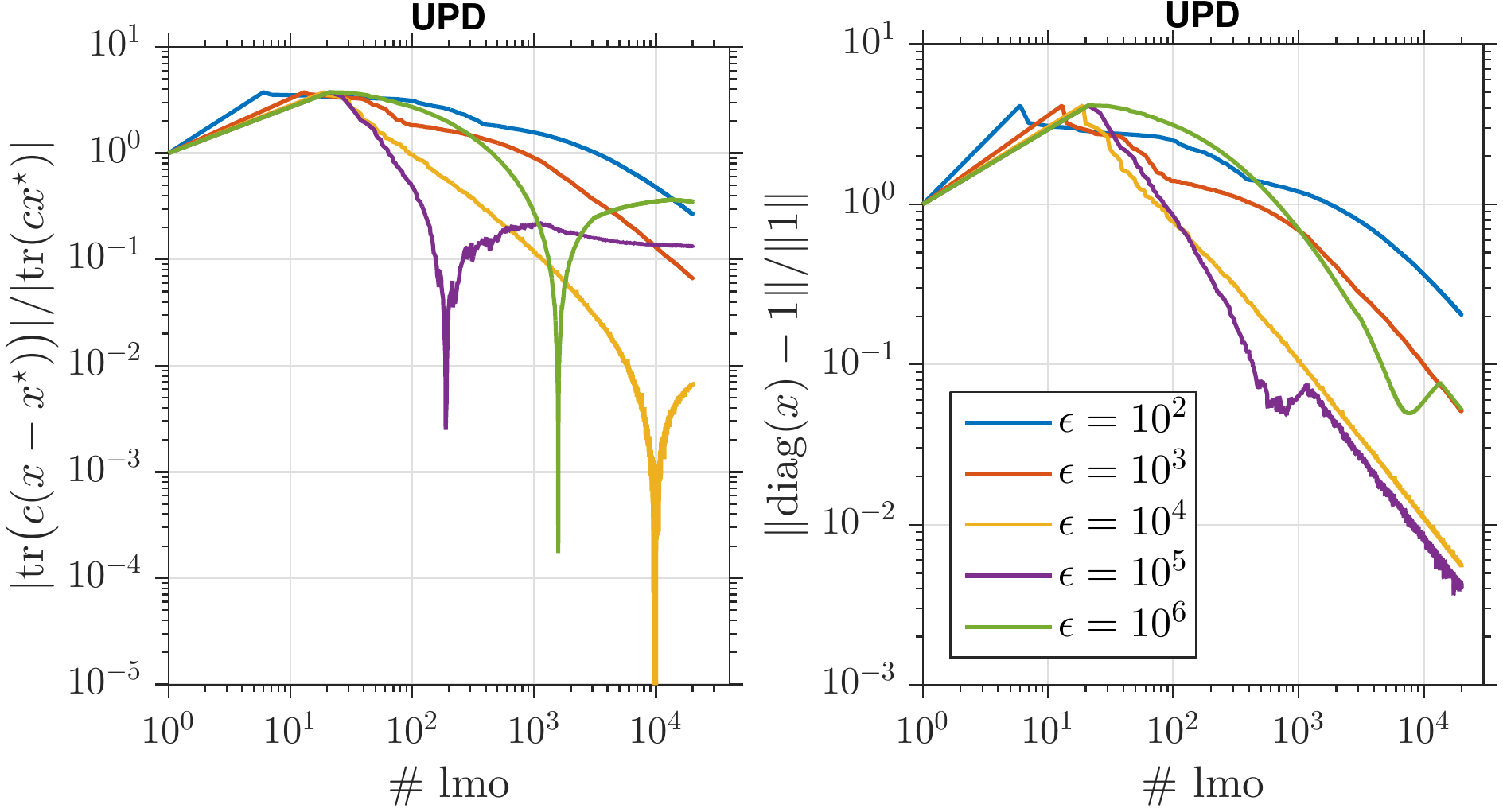}~~~~
\includegraphics[width=0.47\textwidth]{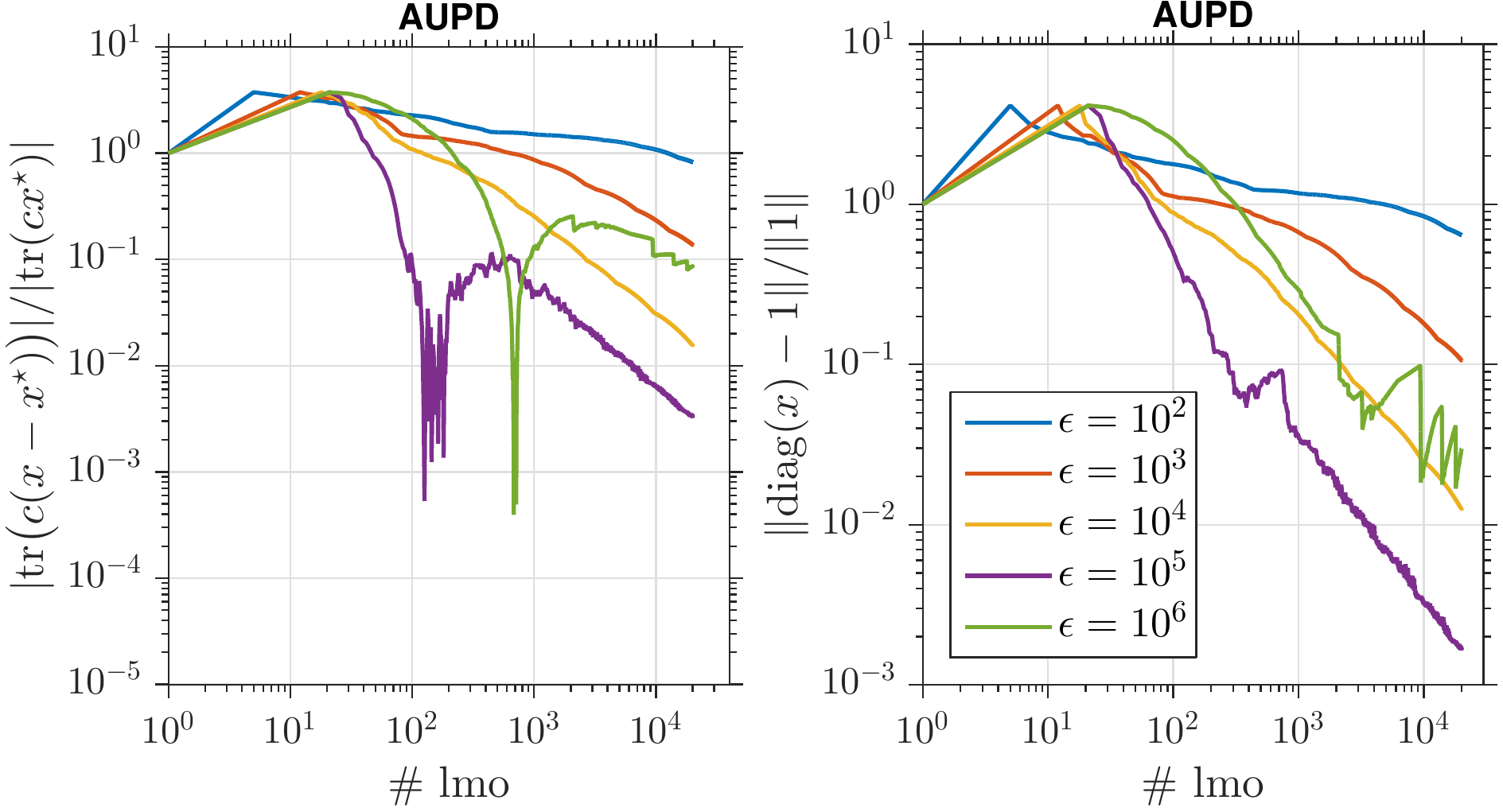}\\[1em]
\includegraphics[width=0.47\textwidth]{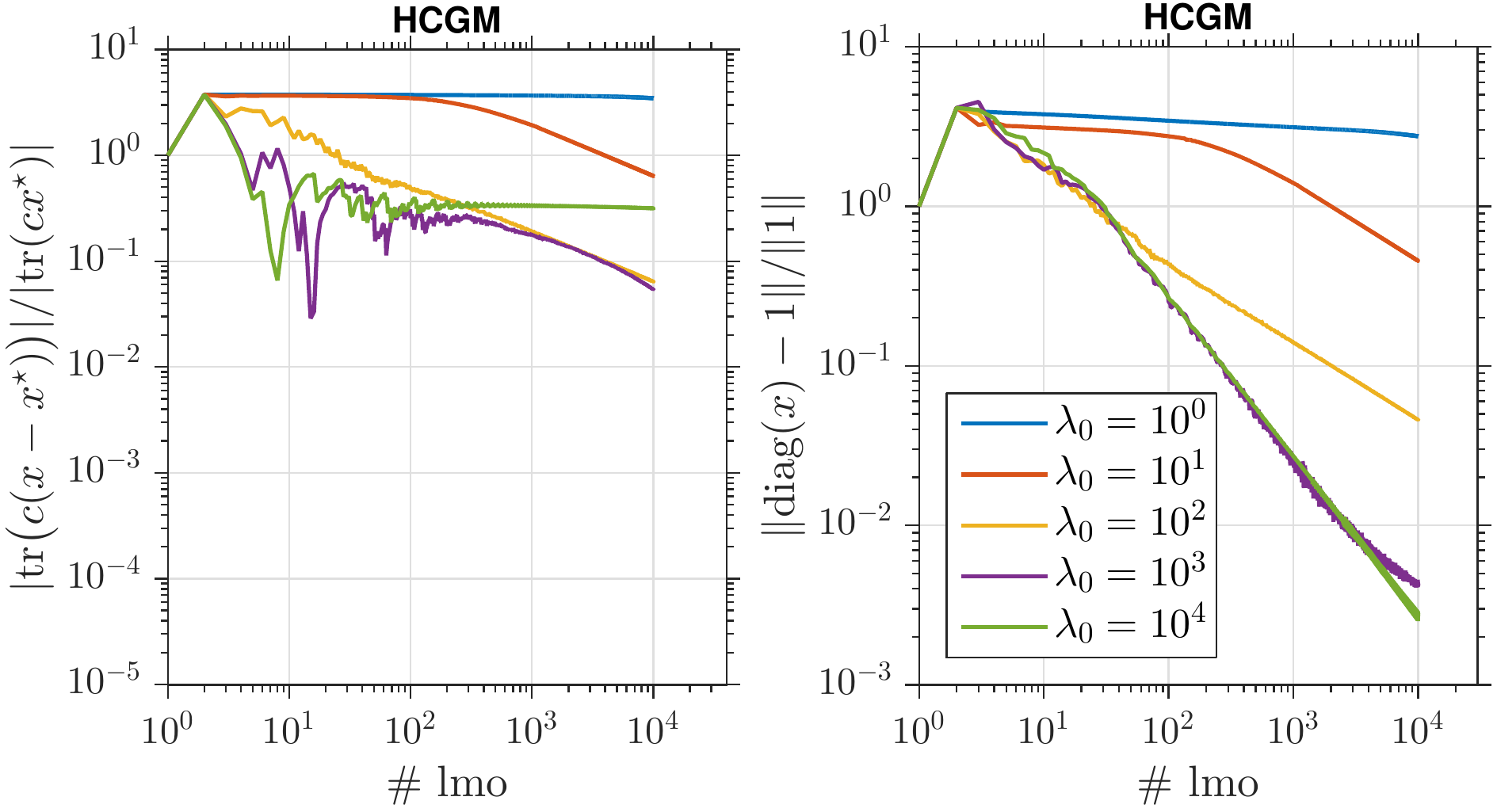}~~~~
\includegraphics[width=0.47\textwidth]{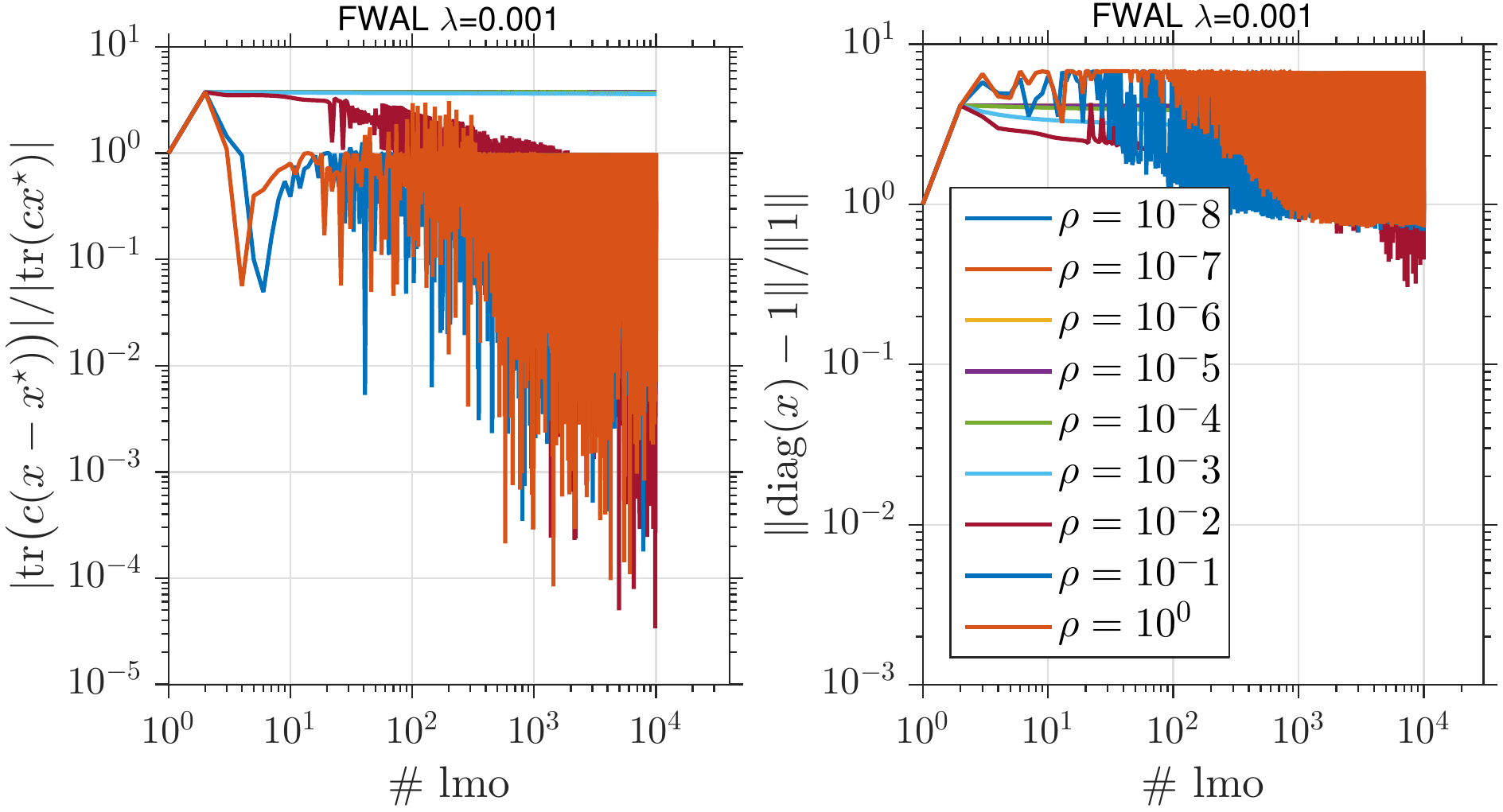}\\[1em]
\includegraphics[width=0.47\textwidth]{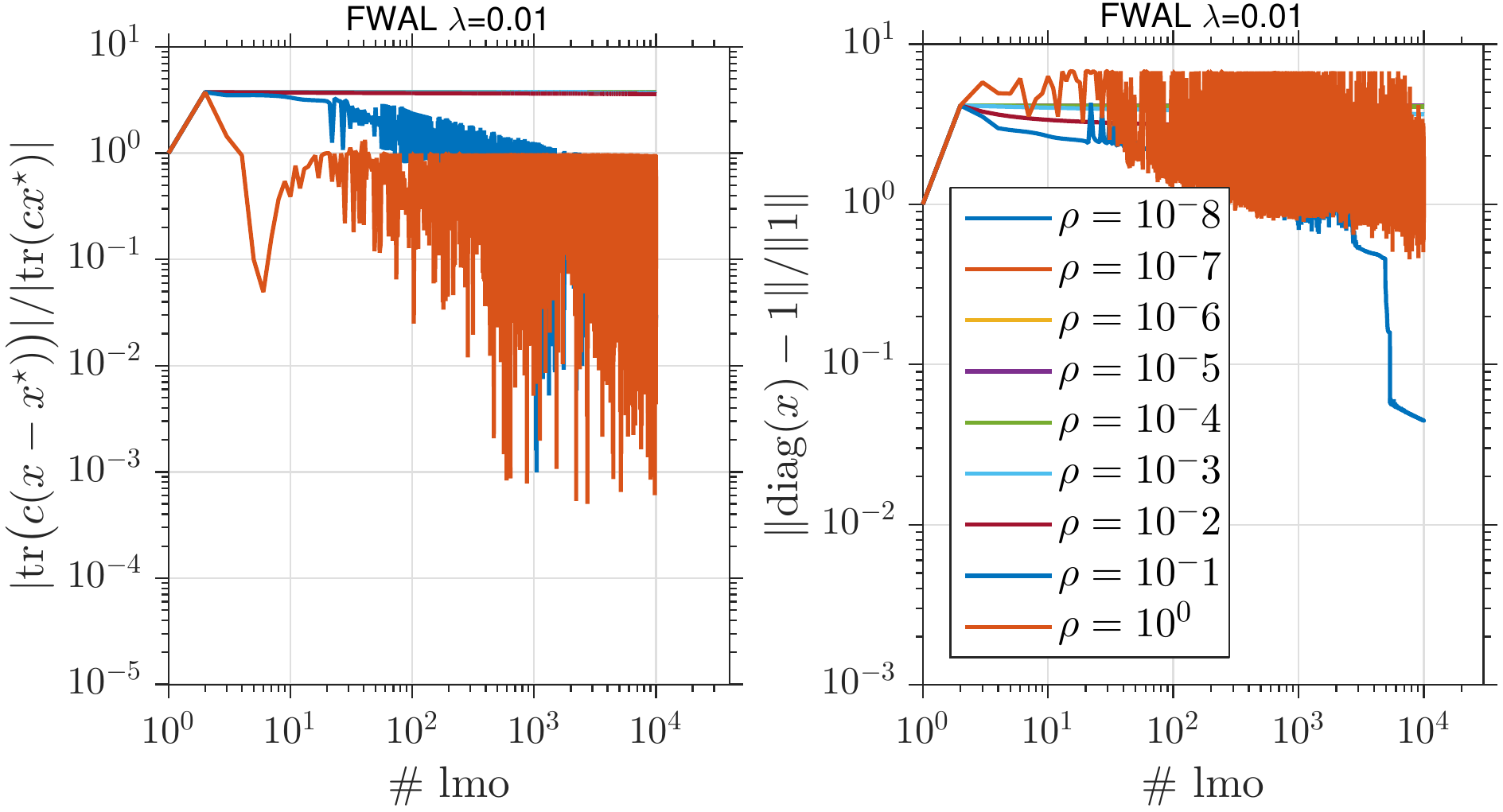}~~~~
\includegraphics[width=0.47\textwidth]{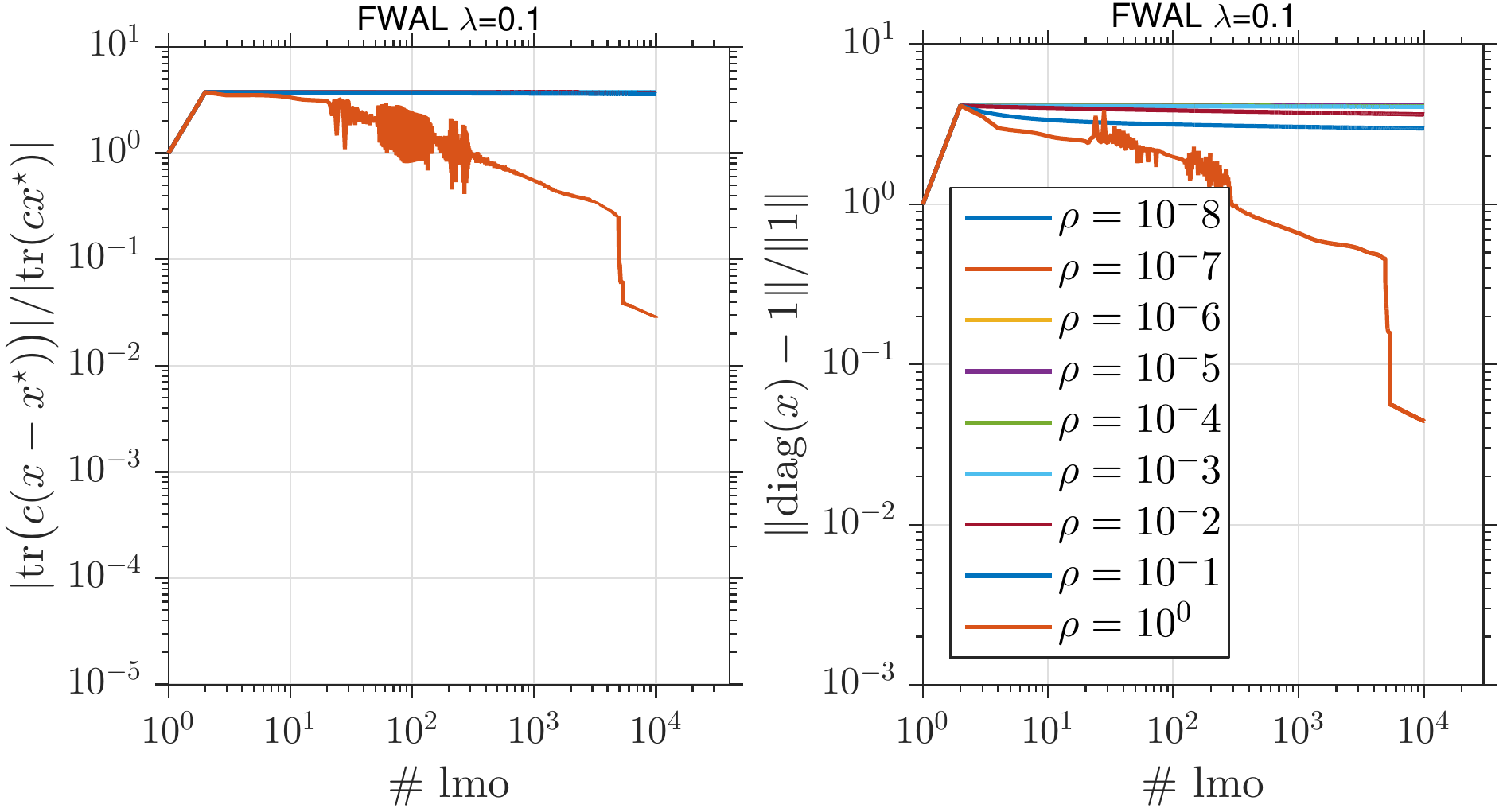}
        \label{fig:maxcut-GD97_b-app1}
\end{figure}

\begin{figure}[ht]
\includegraphics[width=0.47\textwidth]{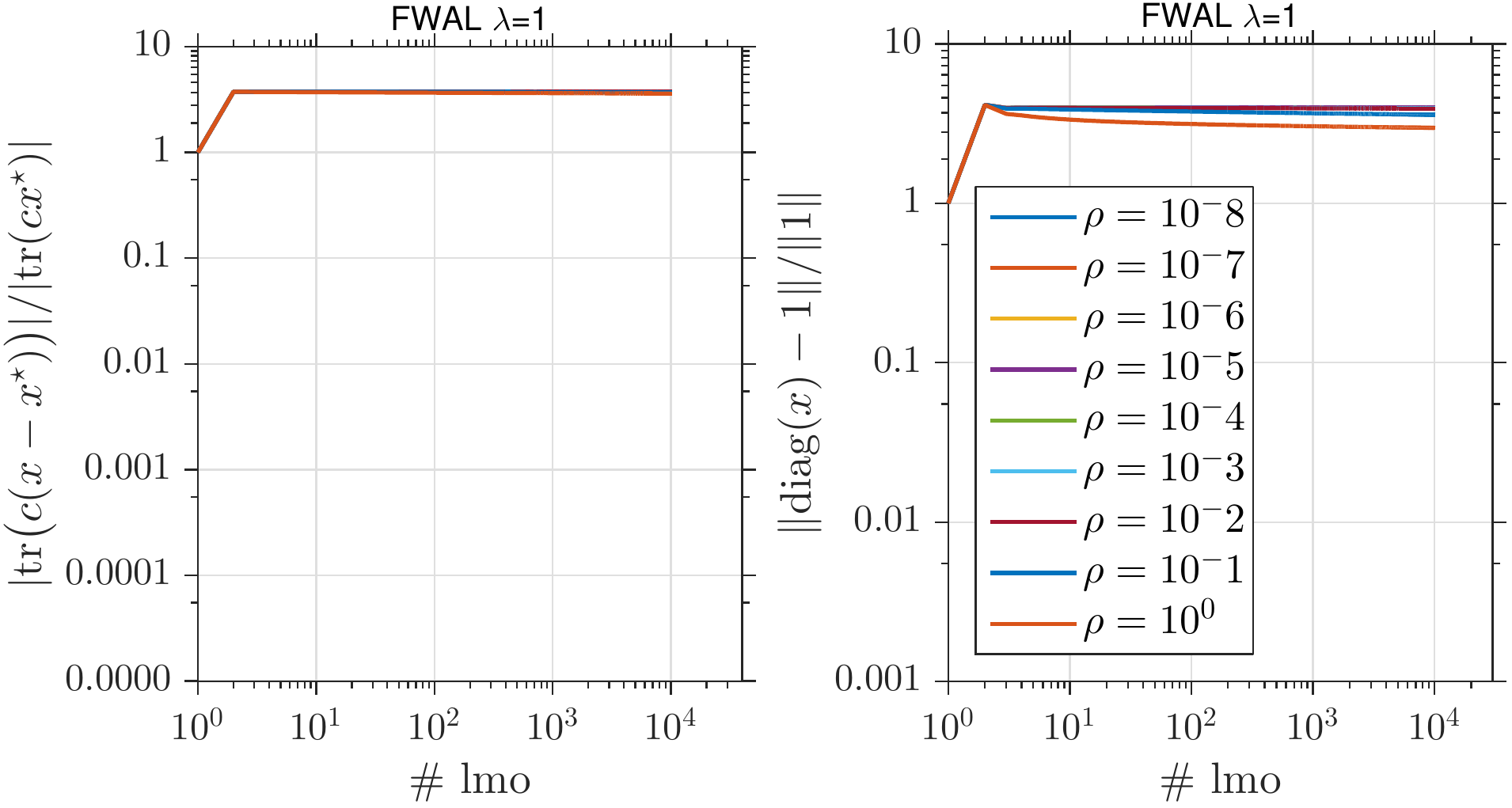}~~~~
\includegraphics[width=0.47\textwidth]{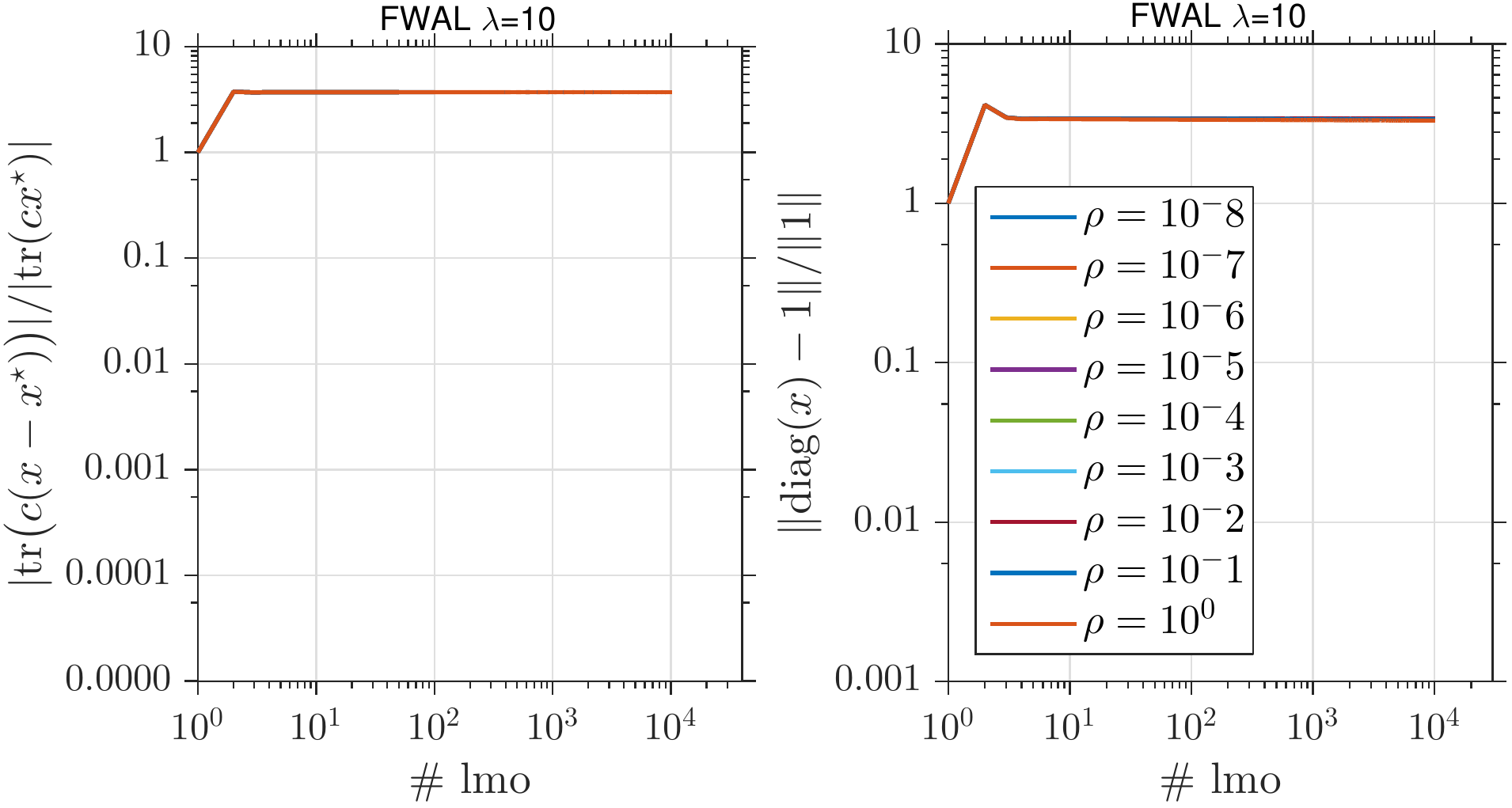}\\[1em]
\includegraphics[width=0.47\textwidth]{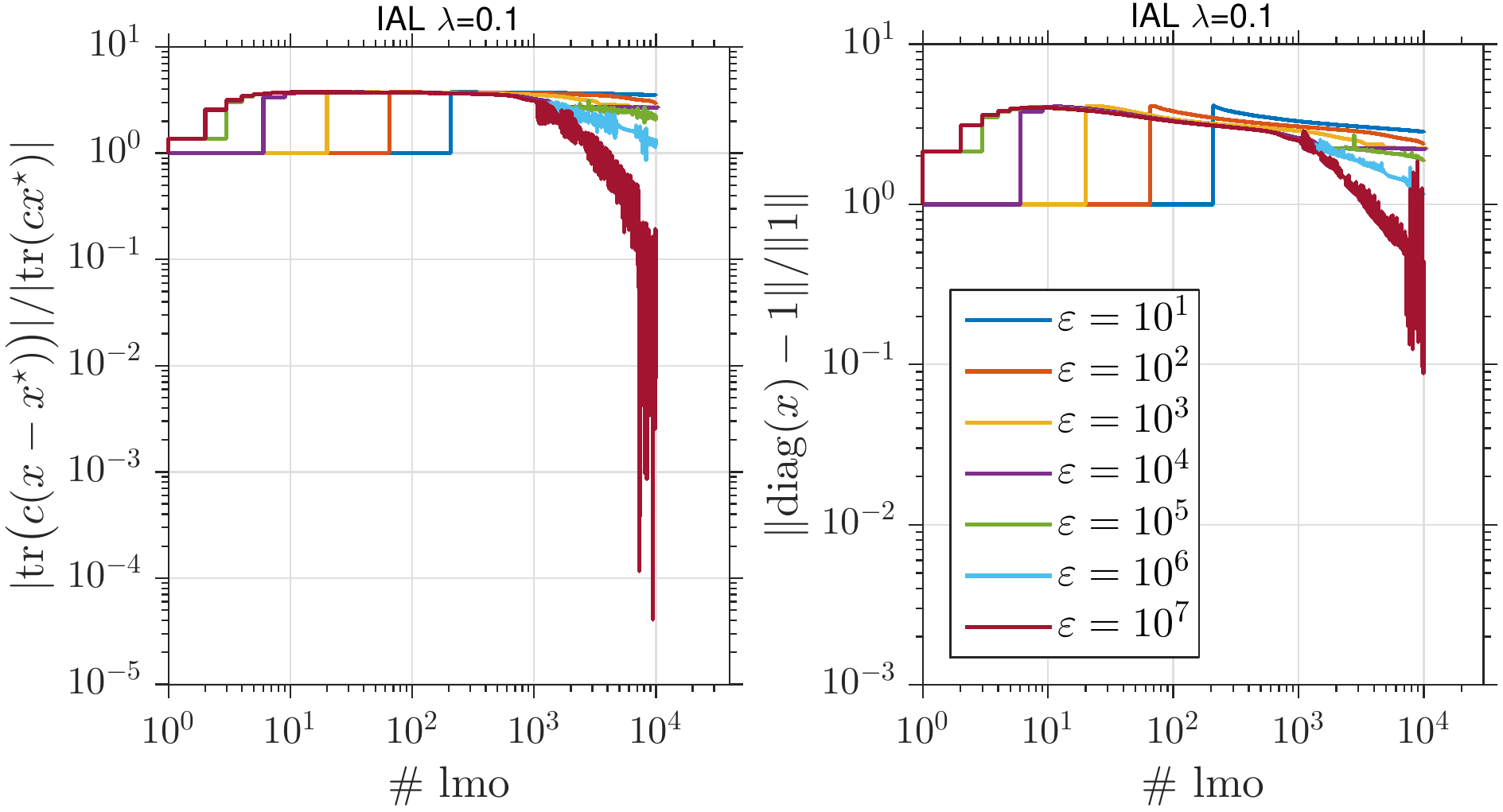}~~~~
\includegraphics[width=0.47\textwidth]{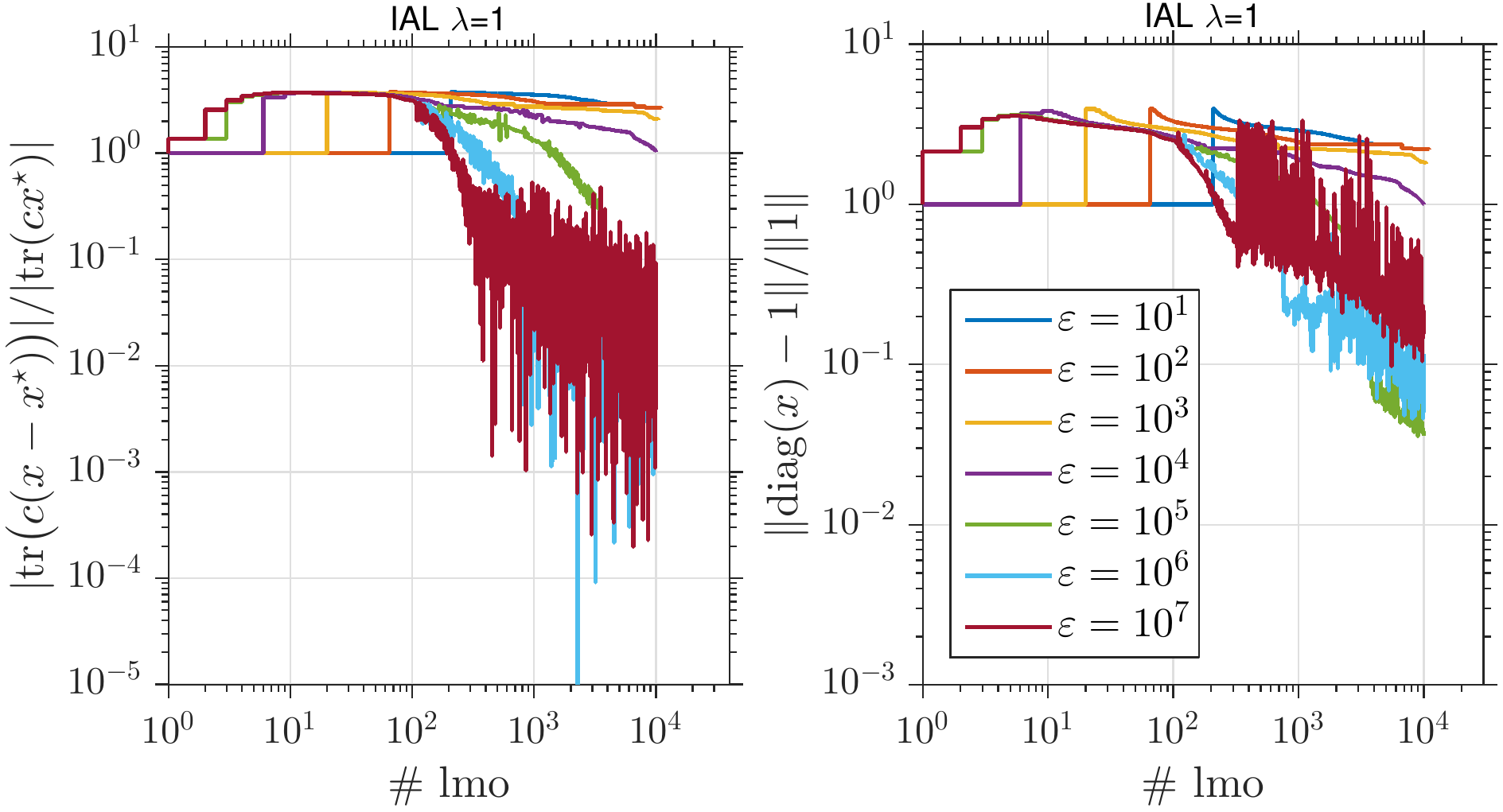}\\[1em]
\includegraphics[width=0.47\textwidth]{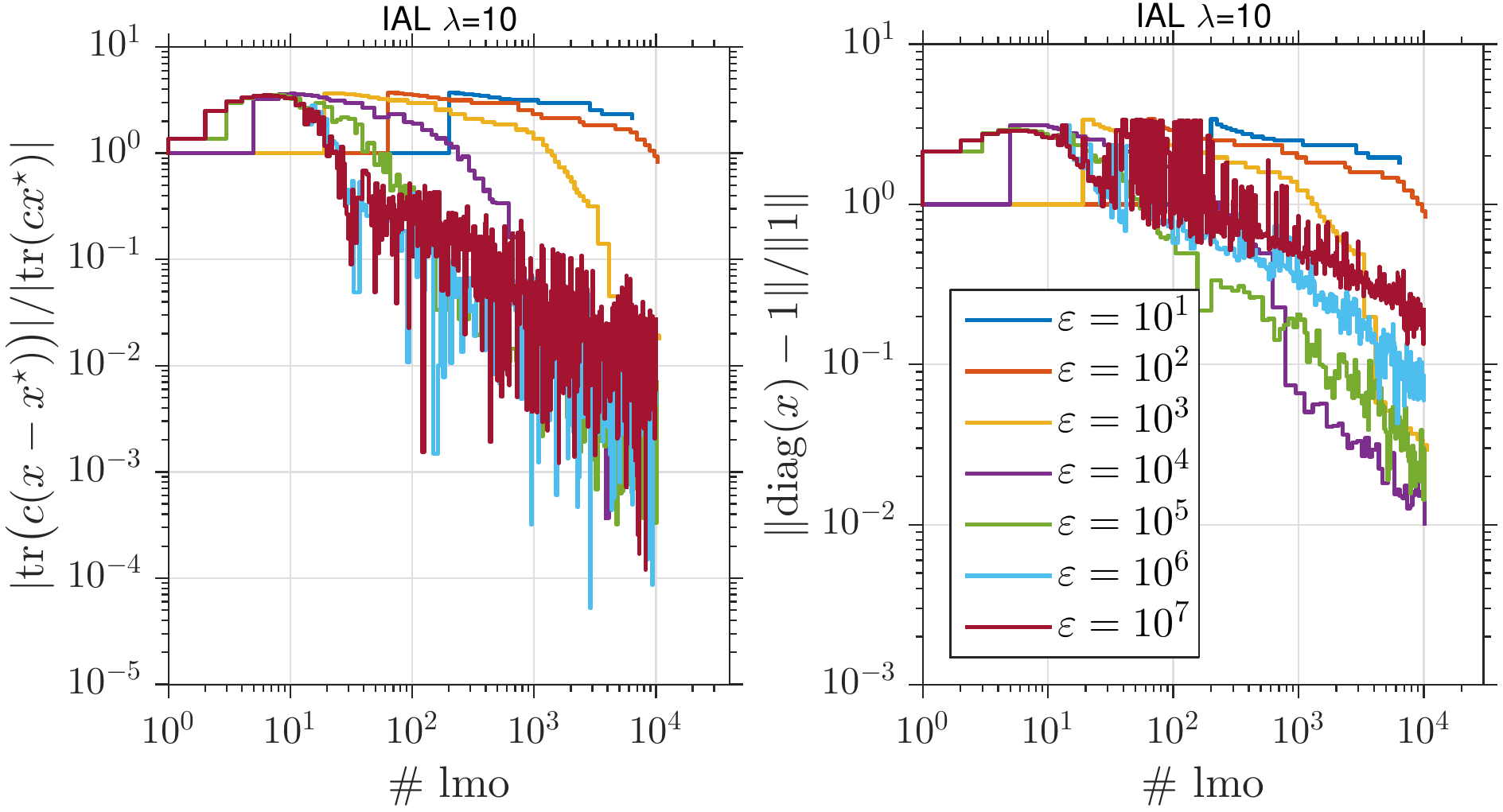}~~~~
\includegraphics[width=0.47\textwidth]{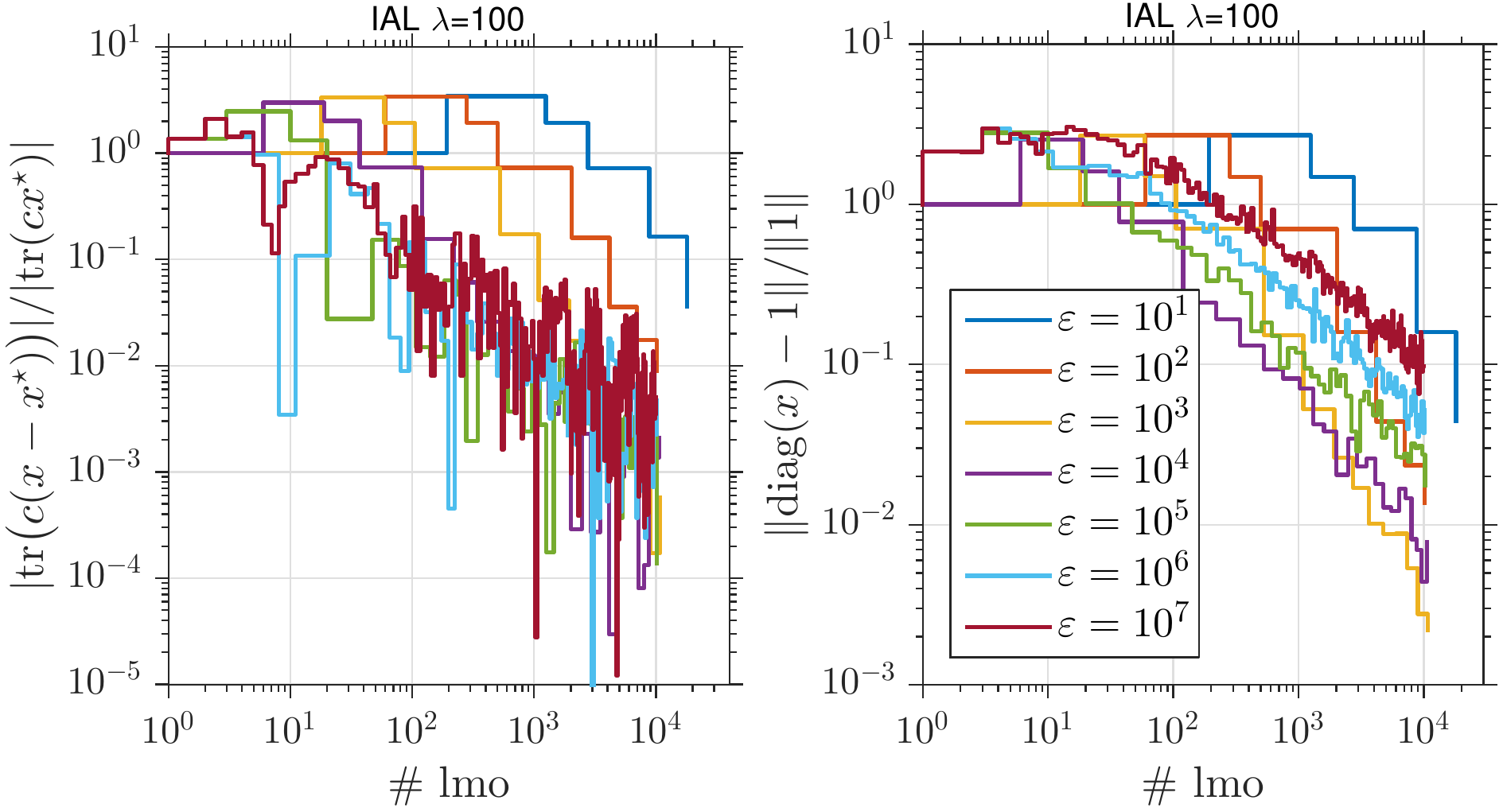}\\[1em]
\includegraphics[width=0.47\textwidth]{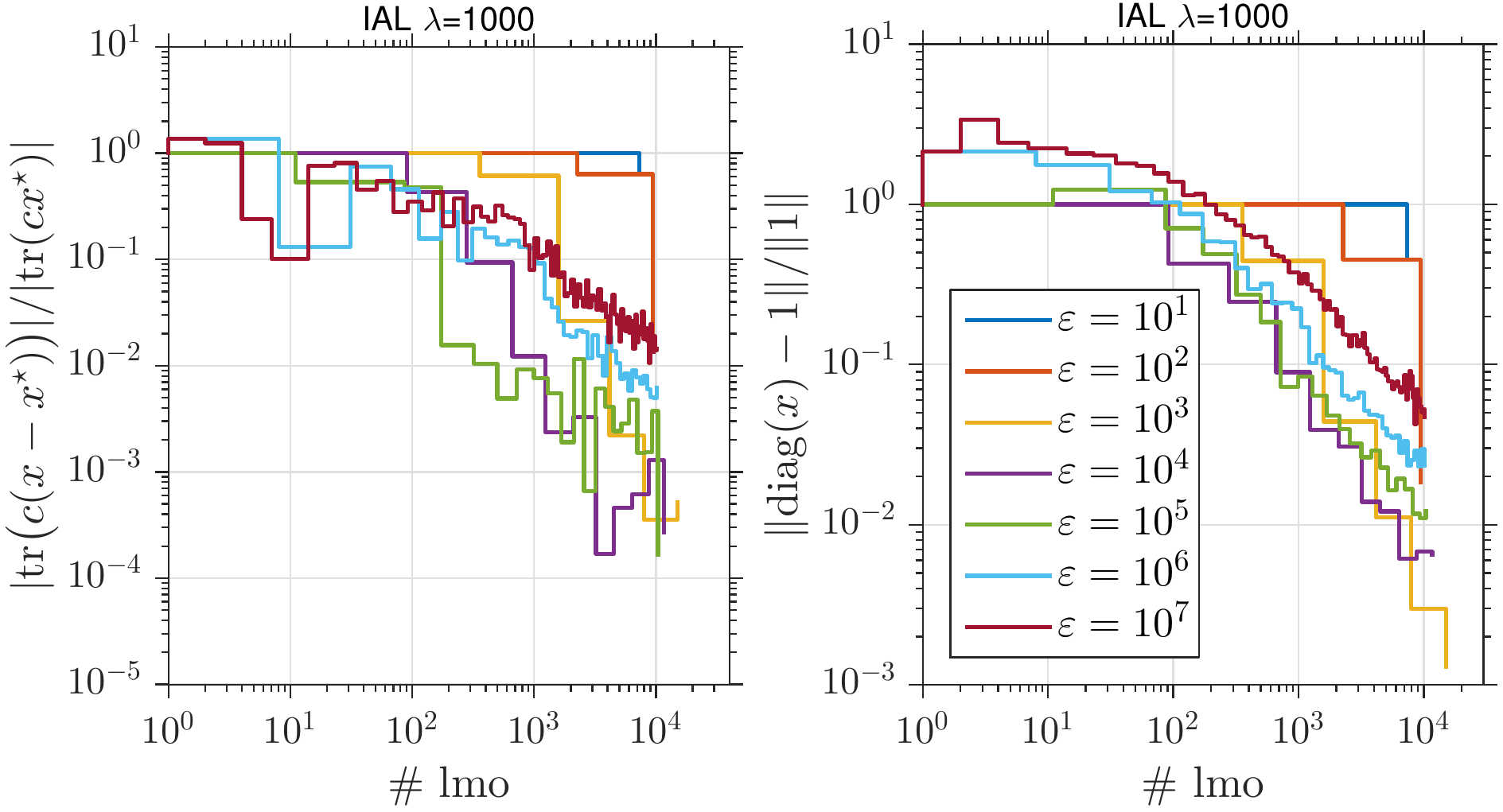}~~~~
\includegraphics[width=0.47\textwidth]{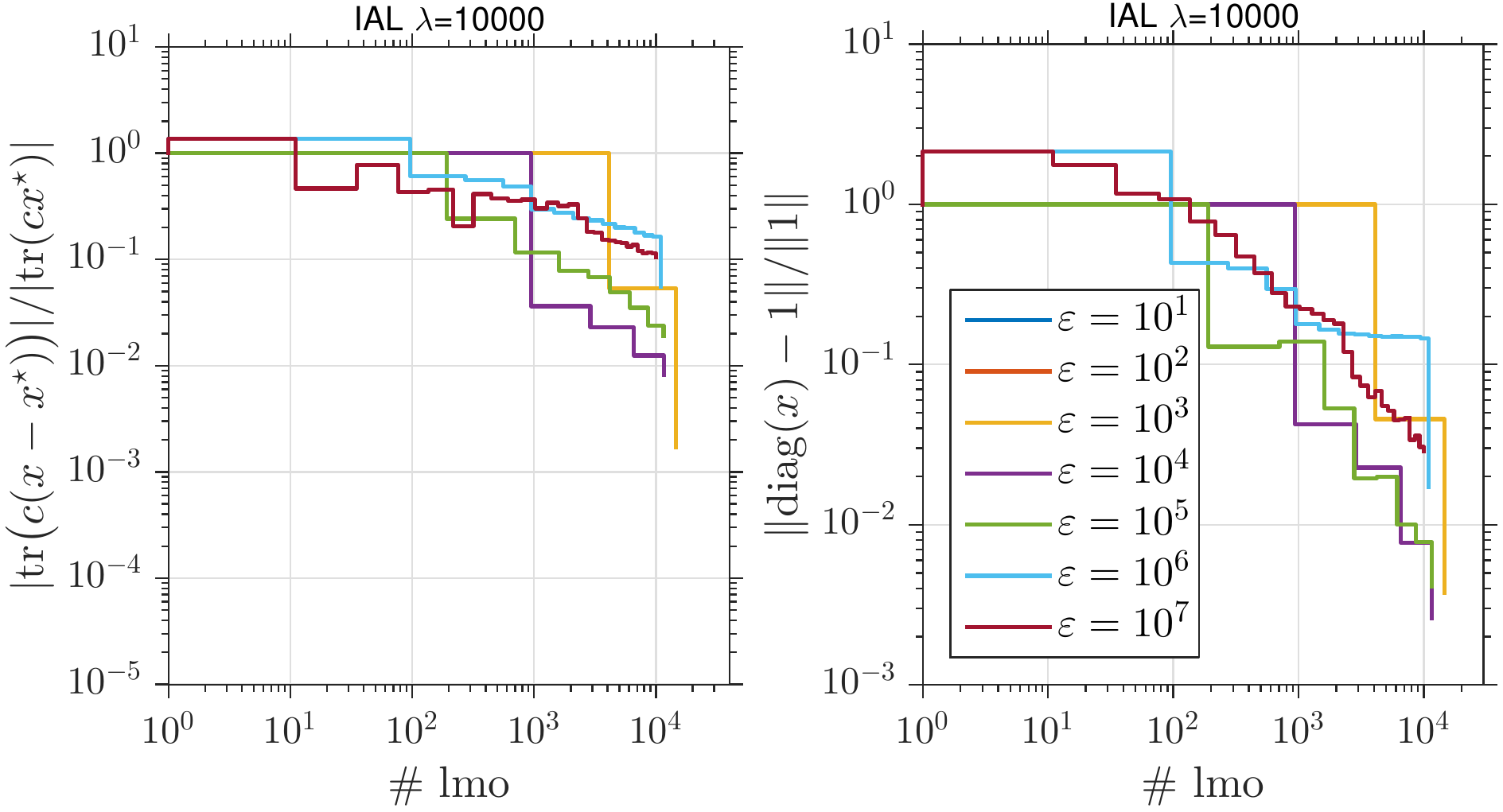}
        \label{fig:maxcut-GD97_b-app2}
\end{figure}

\end{document}